\documentclass[11pt,a4paper]{amsart}
\usepackage[T1]{fontenc}
\usepackage{amssymb,mathtools}

\usepackage[looser]{newpxtext}
\usepackage[smallerops]{newpxmath}

\usepackage{microtype}

\usepackage{enumitem}

\usepackage{geometry}
\usepackage[hidelinks,pagebackref]{hyperref}
\renewcommand*{\backref}[1]{}
\renewcommand*{\backrefalt}[4]{%
\ifcase #1
Not cited%
\or
(Cited on page~#2)%
\else
(Cited on pages~#2)%
\fi
}

\allowdisplaybreaks

\newtheorem{theorem}{Theorem}[section]
\newtheorem{proposition}[theorem]{Proposition}
\newtheorem{lemma}[theorem]{Lemma}
\newtheorem{corollary}[theorem]{Corollary}

\newtheorem{maintheorem}{Theorem}

\theoremstyle{definition}

\newtheorem{example}[theorem]{Example}
\newtheorem{question}[theorem]{Question}
\theoremstyle{remark}
\newtheorem{remark}[theorem]{Remark}

\usepackage{etoolbox}
\makeatletter
\patchcmd{\@settitle}
    {\uppercasenonmath\@title}
    {}
    {}{\typeout{failed}}
\patchcmd{\@setauthors}
    {\MakeUppercase{\authors}}
    {\authors}
    {}{\typeout{failed}}
\makeatother

\newcommand{\C}{\mathbb C}
\newcommand{\T}{\mathbb T}
\newcommand{\cf}{\mathcal F}
\newcommand{\Orb}{\operatorname{Orb}}
\newcommand{\Trans}{\operatorname{Trans}}
\newcommand{\HC}{\operatorname{HC}}
\newcommand{\HCF}{\operatorname{HC}_{\cf}}
\newcommand{\bs}{\boldsymbol}

\title[Weighted orbits of hypercyclic vectors]{Weighted orbits of hypercyclic vectors on F-spaces}

\author[J. Li]{Jian Li}
\address[J. Li]{Institute for  Mathematical Sciences and Artificial Intelligence \& Department of Mathematics,
	Shantou University, Shantou, 515821, Guangdong, China}
\email{lijian09@mail.ustc.edu.cn}
\urladdr{https://orcid.org/0000-0002-8724-3050}

\author[Y. Ruan]{Yonghang Ruan}
\address[Y. Ruan]{Department of Mathematics,
	Shantou University, Shantou, 515821, Guangdong, China}
\email{yonghangruan@outlook.com}
\urladdr{https://orcid.org/0009-0000-2190-6627}

\author[Y. Yang]{Yini Yang}
\address[Y. Yang]{Department of Mathematics,
	Shantou University, Shantou, 515821, Guangdong, China}
\email{ynyangchs@foxmail.com}
\urladdr{https://orcid.org/0000-0001-6564-2213}

\date{\today}



\begin{document}

\begin{abstract}
We study $\cf$-hypercyclicity of weighted orbits $\lambda_nT^nx$ for continuous linear operators $T$ on complex F-spaces, where $\cf$ is a Furstenberg family with the cut--shift--and--paste (CuSP) property.
We characterize the scalar sequences for which weighted hypercyclicity always implies ordinary hypercyclicity: their nonzero values must be bounded and bounded away from zero.
This condition also guarantees weight removal for every Furstenberg family with the CuSP property.
For the converse inclusion, we prove that, for any Furstenberg family $\cf$ with the CuSP property, a bounded scalar sequence preserves all $\cf$-hypercyclic vectors if its backward-shift $\omega$-limit system is distal and does not contain the zero sequence.
For chaotic operators, minimality of this system suffices.
In particular, every complex nilsequence which is not identically zero preserves $\cf$-hypercyclicity, and equality holds when its nonzero values are bounded away from zero.
We also obtain simultaneous weighted return results and stability under unimodular weights with asymptotically trivial multiplicative differences.
\end{abstract}

\maketitle

\section{Introduction}
Let $X$ be a complex F-space, that is, a complex topological vector space admitting a complete translation-invariant metric. 
A map $T\colon X\to X$ is called an \emph{operator} if it is continuous and linear. 
For $x\in X$ and $U\subset X$, put
\[
N_T(x,U)=\{n\in\mathbb N_0\colon T^nx\in U\}.
\]
The vector $x$ is \emph{hypercyclic} for $T$ if $N_T(x,U)$ is nonempty for every nonempty open subset $U$ of $X$.
We denote the set of these vectors by $\HC(T)$. 
We consider the effect of prescribing a scalar at each time of an orbit.
The rotation theorem of Le\'on--Saavedra and M\"uller states that replacing $T$ by $\lambda T$, for $|\lambda|=1$, preserves its hypercyclic vectors \cite{LM04}.
Thus the geometric weights $\lambda^n$ preserve hypercyclicity.
\begin{theorem}[{\cite[Corollary 3]{LM04}}]
Let $T$ be an operator on a complex F-space and let $\lambda\in\mathbb C$ with $|\lambda|=1$.
Then $\HC(T)=\HC(\lambda T)$.
\end{theorem}
Bayart and Costakis extended this assertion to polynomial phases.
\begin{theorem}[{\cite[Theorem 1.1]{BC14}}]\label{thm:bc14}
Let $T$ be an operator on a complex F-space $X$ and let $P$ be a real polynomial.
Then a vector $x\in X$ is hypercyclic for $T$ if and only if
$\{e^{2\pi iP(n)}T^nx\colon n\in\mathbb N_0\}$ is dense in $X$.
\end{theorem}
They also proved the corresponding assertion for frequent hypercyclicity and showed that arbitrary unimodular weights need not preserve hypercyclicity \cite{BC14}.
Consequently, a condition on the moduli alone cannot settle the preservation problem.
The order in which the scalars occur must also be taken into account.

A related problem allows a scalar to be chosen freely from a fixed set $\Gamma\subset\mathbb C$ at each time.
Charpentier, Ernst and Menet characterized the nonempty sets $\Gamma$ for which density of $\{\gamma T^nx\colon \gamma\in\Gamma,\ n\in\mathbb N_0\}$ is equivalent to hypercyclicity of $x$ for all operators on infinite-dimensional complex Banach spaces: $\Gamma\setminus\{0\}$ must be nonempty, bounded and bounded away from zero \cite{CEM16}.
They noted that their result also holds on Fr\'echet spaces.
Alves, Botelho and F\'avaro subsequently studied this scalar-set problem with return times belonging to a Furstenberg family \cite{ABF26}.

The preservation direction leads instead to a compact dynamical system associated with the scalar sequence.
Compact models have been used in linear dynamics by Dayan and Glasner \cite{DG15}, and related Polish-dynamical methods appear in the work of Grivaux and Matheron on invariant measures \cite{GM14}.
Disjointness from minimal distal compact dynamical systems will supply the joint orbit density needed to synchronize the prescribed scalars with the operator orbit.

In this paper, we investigate hypercyclic vectors with respect to Furstenberg families of subsets of $\mathbb N_0$.
Recall that a \emph{(Furstenberg) family} is a nonempty collection $\cf$ of nonempty subsets of $\mathbb N_0$
that is upward closed under inclusion.
Let $T$ be an operator on a complex F-space $X$. 
We denote 
\[\HC_{\cf}(T):=\left\{x\in X\colon N_T(x,U)\in\cf, \forall\text{ nonempty open }U\subset X\right\}\]
and call the elements of $\HC_{\cf}(T)$ the \emph{$\cf$-hypercyclic} vectors of $T$.
The families of infinite sets and of sets with positive lower, upper and upper Banach density give ordinary, frequent, upper-frequent and reiterative hypercyclicity, respectively.
The density conventions and the combinatorial property of families used below are recalled in Section 2. 

For a complex sequence $\bs\lambda=(\lambda_n)_n$, define
\[
N_T^{\bs\lambda}(x,U)=\{n\in\mathbb N_0\colon \lambda_nT^nx\in U\},
\]
and let 
\[
    \HCF^{\bs\lambda}(T): =\left\{x\in X\colon N_T^{\bs\lambda}(x,U)\in\cf,\ \forall\text{ nonempty open }U \subset X\right\}.
\]
For the family $\cf_{\mathrm{inf}}$ of infinite sets, write $\HC^{\bs\lambda}(T)$.
This means precisely that the weighted orbit $\{\lambda_nT^nx\colon n\in\mathbb N_0\}$ is dense.

The question is therefore the following.
\begin{question}\label{main-question}
Let $\cf$ be a family. 
For which scalar sequences $\bs\lambda$ does
\[
\HCF(T)=\HCF^{\bs\lambda}(T)
\]
hold for every complex F-space $X$ and every operator $T\colon X\to X$? 
\end{question}
Write $\bs0=(0,0,\dotsc)$ for the scalar zero sequence.
Our two main results treat the inclusions separately.
The removal of weights depends only on their nonzero values, whereas preservation is governed by their shift dynamics.

\begin{maintheorem}\label{thm:removal}
Let $\bs\lambda$ be a complex sequence.
Then the following assertions are equivalent.
\begin{enumerate}
\item $\HC^{\bs\lambda}(T)\subset\HC(T)$ for every complex F-space $X$ and every operator $T\colon X\to X$;
\item the nonzero values of $\bs\lambda$ are bounded and bounded away from zero. 
\end{enumerate}
Moreover, under the above equivalent conditions, one has $\HCF^{\bs\lambda}(T)\subset\HCF(T)$ for every complex F-space $X$, every operator $T\colon X\to X$ and every family $\cf$ with the CuSP property.
\end{maintheorem}

The CuSP property is defined in Section 2 and includes the standard density families mentioned above.
For (1)$\Rightarrow$(2) in Theorem \ref{thm:removal}, we construct a weighted hypercyclic vector for the backward shift $B$ on $\mathbb C^{\mathbb N_0}$ whose unweighted orbit is not dense whenever either bound on the nonzero weights fails.
For (2)$\Rightarrow$(1), we use our extension of the scalar-set theorem of \cite{CEM16} to complex F-spaces. 
For the assertion of families with the CuSP property, we use the CuSP property to transfer return times through a finite cover of the scalar set.

For a fixed family, the removal problem admits a finer characterization. Building on the proof of Theorem \ref{thm:removal}, we characterize weight removal in purely combinatorial terms through the dual family under additional assumptions (Theorem \ref{thm:F-removal}). 
This yields necessary and sufficient conditions for weight removal in reiterative, upper-frequent, and piecewise-syndetic hypercyclicity (Corollary \ref{coro:F-removal}).

To state the preservation result, let $B$ be the \emph{backward shift} on $\mathbb C^{\mathbb N_0}$ (endowed with the product topology), that is, $B(z_0,z_1,\dotsc)=(z_1,z_2,\dotsc)$.
For a bounded complex sequence $\bs\lambda$, its \emph{$\omega$-limit set} is
\[
\Omega_{\bs\lambda}: =\bigcap_{N\in\mathbb N_0}
\overline{\{B^n\bs\lambda\colon n\geq N\}}.
\]
This is a nonempty compact subset of $\mathbb C^{\mathbb N_0}$ with $B\Omega_{\bs\lambda}=\Omega_{\bs\lambda}$. 

\begin{maintheorem}\label{thm:shift-criterion}
Let $\bs\lambda$ be a bounded complex sequence.
If $(\Omega_{\bs\lambda},B)$ is distal with $\bs 0\notin \Omega_{\bs\lambda}$, then for every complex F-space $X$, every operator $T\colon X\to X$ and every family $\cf$ with the CuSP property, one has
\[
      \HCF(T)\subset \HCF^{\bs\lambda}(T).
\]
\end{maintheorem}
Theorem \ref{thm:shift-criterion} is a sufficient criterion, which does not characterize all sequences which preserve hypercyclicity.
Together with Theorem \ref{thm:removal}, 
it gives equality $\HCF(T)=\HCF^{\bs\lambda}(T)$ whenever the nonzero weights lie in a fixed annulus and the stated condition on the $\omega$-limit system holds.

The proof of Theorem \ref{thm:shift-criterion} follows a different approach.
Shkarin's result \cite{S08} implies that hypercyclic operators are weakly scattering.
The disjointness result of Huang and Ye \cite{HY12} then gives joint orbit density with each minimal distal subsystem of the scalar $\omega$-limit set.
Compactness reduces these pointwise density statements to finitely many time shifts, and the CuSP property puts their return sets back together.
For chaotic operators, disjointness from every minimal system allows minimality to replace distality in the scalar hypothesis.
We also give versions for direct sums and for simultaneous returns of a single vector under several weights.

Nilsequences provide a class of scalar weights to which this argument applies.
Every nilsequence can be realized as a continuous observation of a distal compact dynamical system. 
It follows that every nilsequence which is not identically zero preserves $\cf$-hypercyclicity.
Such a sequence may have zero terms or nonzero terms tending arbitrarily close to zero. 
The additional lower bound is needed only for the reverse inclusion.
We finally combine weights with distal $\omega$-limit systems with unimodular sequences whose iterated multiplicative differences tend to one.

The paper is organized as follows.
In Section 2, we introduce the required facts about families, compactifications and disjointness.
In Section 3, we prove Theorem \ref{thm:removal} and establish further combinatorial characterizations of weight removal.
In Section 4, we prove Theorem \ref{thm:shift-criterion} and its extensions. 
In Section 5, to give applications of Theorems \ref{thm:removal} and \ref{thm:shift-criterion}, we treat periodic weights, polynomial phases, nilsequences and unimodular sequences whose iterated multiplicative differences tend to one. 

\section{Furstenberg families and Polish dynamical systems}

In this section, we collect some basic facts about
Furstenberg families, Polish dynamical systems, and weighted hypercyclicity that will be needed in the sequel.
\subsection{Furstenberg families}
Write $\mathbb N=\{1,2,\dotsc\}$ and $\mathbb N_0=\mathbb N\cup\{0\}$.
A nonempty collection $\cf$ of nonempty subsets of $\mathbb N_0$ is called a \emph{Furstenberg family} (or \emph{family} for short) if it is hereditary upward: for any $A\subset B\subset\mathbb N_0$ with $A\in\cf$, one has $B\in\cf$. 
The \emph{dual family} of $\cf$ is given by 
\[
    \cf^*=\left\{F\subset\mathbb N_0\colon F\cap A\neq\emptyset,\ \forall A\in\cf\right\}. 
\]
Detailed discussions of Furstenberg families can be found in \cite{F81}.

We recall some standard families of subsets of $\mathbb N_0$.
Let $\mathcal F_{\mathrm{inf}}$ denote the family of all infinite subsets of
$\mathbb N_0$.
A set $A\subset\mathbb N_0$ is \emph{cofinite} if $\mathbb N_0\setminus A$ is finite.
We denote the family of all cofinite subsets of $\mathbb N_0$ by $\cf_{\mathrm{cf}}$.
A set $F\subset\mathbb N_0$ is \emph{syndetic} if it has bounded gaps, and
\emph{thick} if it contains intervals of arbitrary length.
We denote by $\mathcal F_{\mathrm s}$ and $\mathcal F_{\mathrm t}$
the families of syndetic and thick subsets of $\mathbb N_0$, respectively.
A set is \emph{piecewise syndetic} if it is the intersection of a syndetic
set and a thick set, and we denote the family of all such sets by
$\mathcal F_{\mathrm{ps}}$. 
A set $F\subset\mathbb N_0$ is \emph{thickly syndetic} if, for every $L\in\mathbb N$, the set
$\{n\in\mathbb N_0:[n,n+L-1]\cap\mathbb N_0\subset F\}$ is syndetic.
We denote the family of all thickly syndetic subsets of $\mathbb N_0$ by $\mathcal F_{\mathrm{ts}}$.

For $A\subset\mathbb N_0$, the \emph{lower density} and \emph{upper density} of $A$ are defined as
\[
\underline d(A)=\liminf_{N\to\infty}\frac{|A\cap[0,N)|}{N},\quad\text{ and }\quad \overline d(A)=\limsup_{N\to\infty}\frac{|A\cap[0,N)|}{N}, 
\]
and the \emph{lower Banach density} and \emph{upper Banach density} of $A$ are defined as
\[
\underline{\operatorname{Bd}}(A) =\lim_{N\to\infty}\inf_{m\in\mathbb{N}_0}
\frac{|A\cap[m,m+N)|}{N},\quad\text{ and }\quad \overline{\operatorname{Bd}}(A) =\lim_{N\to\infty}\sup_{m\in\mathbb{N}_0}
\frac{|A\cap[m,m+N)|}{N}.
\]
Intervals in these formulas are intervals of integers.
We denote the families of sets with positive lower density, positive upper density, and positive upper Banach density by $\cf_{\mathrm{pld}}$, $\cf_{\mathrm{pud}}$, and $\cf_{\mathrm{pubd}}$, respectively.
The families of sets with lower density one and lower Banach density one are denoted by $\cf_{\mathrm{ld1}}$ and $\cf_{\mathrm{lbd1}}$, respectively.
It is known that $\cf_{\mathrm{inf}}^*=\cf_{\mathrm{cf}}$, $\cf_{\mathrm{pubd}}^*=\cf_{\mathrm{lbd1}}$, $\cf_{\mathrm{pud}}^*=\cf_{\mathrm{ld1}}$ and $\cf_{\mathrm{ps}}^*=\cf_{\mathrm{ts}}$. 

We say that a family $\cf$ has the \emph{Ramsey property} if for any $A_1,A_2,\dotsc,A_q\subset\mathbb N_0$ with $A_1\cup A_2\cup\dotsb \cup A_q\in\cf$, one has $A_j\in\cf$ for some $j\in\{1,2,\dotsc, q\}$, or equivalently, we have $F_1\cap F_2\in\cf$ whenever $F_1\in\cf$ and $F_2\in\cf^*$. 
The families $\cf_{\mathrm{inf}}$, $\cf_{\mathrm{pubd}}$ and $\cf_{\mathrm{pud}}$ have the Ramsey property.
The same holds for $\cf_{\mathrm{ps}}$ (see \cite[Theorem 1.24]{F81}).

Following \cite{BGLP22}, we recall the notion of the CuSP property.
A family $\cf$ has the \emph{cut--shift--and--paste property} (or \emph{CuSP property} for short) if for every finite cover $\mathbb N_0=I_1\cup\cdots\cup I_q$, every $r_1,\dotsc,r_q\in\mathbb N_0$, and every $ A\in\cf$,
\[
   \bigcup_{j=1}^q\bigl(r_j+(A\cap I_j)\bigr)\in\cf.
\]
The families $\cf_{\mathrm{inf}},\cf_{\mathrm s},\cf_{\mathrm{pld}},\cf_{\mathrm{pud}}$ and $\cf_{\mathrm{pubd}}$ have this property \cite[Lemma 4.1]{BGLP22}.
We shall also use its equivalent formulation for a finite cover $A=A_1\cup\cdots\cup A_q$ of a member $A\in\cf$:
\[
\bigcup_{j=1}^q(r_j+A_j)\in\cf.
\]
Indeed, extend the cover to $\mathbb N_0$ by adding $\mathbb N_0\setminus A$ to its first member.
Taking $q=1$ shows that the CuSP property implies right-invariance: $A\in\cf$ implies $r+A\in\cf$ for every $r\in\mathbb N_0$.

We recall the notion of block family.
For a family $\cf$, its \emph{block family} $b\cf$ consists of all $S\subset\mathbb N_0$ for which there exists $F\in\cf$ such that every finite set $E\subset F$ has a translate contained in $S$, i.e.\ there exists $t\in\mathbb N_0$ with $t+E\subset S$.
We next show that the block family of a family with the CuSP property also has the CuSP property.

\begin{proposition}
If a family $\cf$ has the CuSP property, then so does $b\cf$.
In particular, $\cf_{\mathrm{ps}}$ has the CuSP property.
\end{proposition}

\begin{proof}
Let $S\in b\cf$.
Choose $F\in\cf$ such that every finite subset of $F$ has a translate contained in $S$.
It suffices to treat a finite partition $S=S_1\cup\cdots\cup S_q$ with assigned shifts $r_1,r_2,\dotsc,r_q\in\mathbb N_0$.

We first find a map $c\colon F\to\{1,\dotsc,q\}$ such that, for every finite subset $E$ of $F$, there exists $t\in\mathbb N_0$ with $t+f\in S_{c(f)}$ for every $f\in E$.
It suffices to consider the case where $F$ is infinite and write $F=\{f_1,f_2,\dotsc\}$.
For each $m\in\mathbb N$, choose $t_m\in\mathbb N_0$ such that
$t_m+\{f_1,\dotsc,f_m\}\subset S$.
For $i=1,\dotsc,m$, choose $c_m(f_i)\in\{1,\dotsc,q\}$ such that $t_m+f_i\in S_{c_m(f_i)}$.
Since there are only finitely many colors, by successively passing to subsequences, we can choose $c(f_i)\in\{1,\dotsc,q\}$, $i\in\mathbb N$, such that, for every $m\in\mathbb N$, there exists $t\in\mathbb N_0$ satisfying
$t+f_i\in S_{c(f_i)}$ for $i=1,\dotsc,m$.
This gives the required map $c$.
The case where $F$ is finite is immediate.
For $j=1,\dotsc,q$, put $Q_j=\{f\in F\colon c(f)=j\}$.
Then $F=Q_1\cup\cdots\cup Q_q$.
Since $\cf$ has the CuSP property, we have
\[
W: =\bigcup_{j=1}^q(r_j+Q_j)\in\cf.
\]

Let $W_0$ be a finite subset of $W$.
For each $h\in W_0$, there exist $j_h\in\{1,\dotsc,q\}$ and $f_h\in Q_{j_h}$ such that $h=r_{j_h}+f_h$.
Since $f_h\in Q_{j_h}$, we have $c(f_h)=j_h$.
By the choice of $c$, there exists $t\in\mathbb N_0$ such that $t+f_h\in S_{j_h}$ for every $h\in W_0$.
Hence $t+h=r_{j_h}+(t+f_h)\in r_{j_h}+S_{j_h}$ for every $h\in W_0$. Therefore $t+W_0\subset\bigcup_{j=1}^q(r_j+S_j)$.
Thus $\bigcup_{j=1}^q(r_j+S_j)\in b\cf$, and hence $b\cf$ has the CuSP property.
Finally, since $b\cf_{\mathrm s}=\cf_{\mathrm{ps}}$ and $\cf_{\mathrm s}$ has the CuSP property, $\cf_{\mathrm{ps}}$ has the CuSP property.
\end{proof}

\subsection{Polish dynamical systems}

A \emph{Polish dynamical system} is a pair $(X,T)$,
where $X$ is a Polish space and $T\colon X\to X$ is continuous.
Two special cases will be of particular interest: compact dynamical systems, for which $X$ is a compact metrizable space, and linear dynamical systems, for which $X$ is a separable complex F-space and $T$ is an operator on $X$. 

Let $(X,T)$ be a Polish dynamical system. 
A point $x\in X$ is a \emph{transitive point} if 
\[
N_T(x,U):=\{n\in\mathbb N_0\colon T^n x\in U\}\neq\emptyset
\] for every nonempty open subset $U$ of $X$, or equivalently, its orbit 
\[
\Orb(x,T): =\{T^nx\colon n\in\mathbb N_0\} 
\] 
is dense in $X$.
We denote the set of transitive points by $\Trans(T)$.
A Polish dynamical system $(X,T)$ is said to be \emph{(topologically) transitive} if for every pair $U, V$ of nonempty open subsets of $X$, there exists $n\in\mathbb N_0$ such that $T^nU\cap V\neq\emptyset$.
We shall use the following classical form of the Birkhoff transitivity theorem.

\begin{theorem}[see e.g.\ {\cite[Theorem 1.16]{GEP11}}]\label{thm:birkhoff}
    Let $(X,T)$ be a Polish dynamical system.
    Consider the following assertions:
    \begin{itemize}
        \item [(1)] $T$ is transitive;
        \item [(2)]$\Trans(T)$ is a dense $G_\delta$ subset of $X$;
        \item [(3)]$\Trans(T)\neq\emptyset$.
    \end{itemize}
   Then we have (1) $\Leftrightarrow$ (2) $\Rightarrow$ (3).
If $X$ has no isolated points, then (1) $\Leftrightarrow$ (2) $\Leftrightarrow$ (3).
\end{theorem}

To transfer compact dynamical arguments to the Polish setting, we shall use dynamical compactifications.
A \emph{dynamical compactification} (or \emph{compactification} for short) of a Polish dynamical system $(X,T)$ is a compact dynamical system $(\widetilde X,\widetilde T)$ together with an equivariant topological embedding
\[
    J\colon X\hookrightarrow\widetilde X,
    \quad \widetilde T\circ J=J\circ T,
\]
whose image is dense in $\widetilde X$.
It is well known that every Polish dynamical system admits a compactification, see e.g.\ \cite[(7.1)]{O52}.
For completeness, we include the proof here.
\begin{proposition}\label{prop:compactification-exists}
Every Polish dynamical system $(X,T)$ admits a compactification.
\end{proposition}

\begin{proof}
Let $j\colon X\to[0,1]^{\mathbb N_0}$ be a topological embedding into the Hilbert cube.
Put $\mathcal K=\bigl([0,1]^{\mathbb N_0}\bigr)^{\mathbb N_0}$ and let $\sigma$ be the backward shift on $\mathcal K$.
Define $J(x)=\bigl(j(x),j(Tx),j(T^2x),\dotsc\bigr)$.
The zeroth coordinate shows that $J$ is a topological embedding, and $\sigma J(x)=J(Tx)$.
Set $\widetilde X=\overline{J(X)}$ and let $\widetilde T=\sigma|_{\widetilde X}$.
Then $\widetilde X$ is a compact metrizable space such that $J(X)$ is dense.
Thus $(\widetilde X,\widetilde T)$ is a compactification of $(X,T)$.
\end{proof}

We shall identify $X$ with $J(X)$ in a chosen compactification.
Since a Polish subspace of a compact metrizable space is $G_\delta$, $X$ is a dense $G_\delta$ subset of $\widetilde X$.
We shall also use without further mention the following elementary facts: 
if $(\widetilde X,\widetilde T)$ is a compactification of $(X,T)$, then $(\widetilde X,\widetilde T^p)$ is a compactification of $(X,T^p)$ for every $p\in\mathbb N$; 
and if $(\widetilde X,\widetilde T)$ and $(\widetilde Y,\widetilde S)$ are compactifications of $(X,T)$ and $(Y,S)$, respectively, then $(\widetilde X\times\widetilde Y,\widetilde T\times\widetilde S)$ is a compactification of $(X\times Y,T\times S)$.
Transitivity is preserved under taking compactification.
Precisely, we have the following observation.
\begin{proposition}\label{prop:transitive-compactification}
Let $(X,T)$ be a Polish dynamical system with a compactification $(\widetilde X,\widetilde T)$.
Then
\[
    \Trans(T)=\Trans(\widetilde T)\cap X.
\]
In particular, $T$ is transitive if and only if so is $\widetilde{T}$.
\end{proposition}

\begin{proof}
The equality follows from the easy fact that for every $x\in X$ and every nonempty open subset $\widetilde{U}$ of $\widetilde{X}$, we have $N_{\widetilde{T}}(x,\widetilde{U})=N_T(x,\widetilde{U}\cap X)$.
The equivalence of the transitivity of $T$ and $\widetilde{T}$ follows immediately from the equality and Theorem \ref{thm:birkhoff}, since $X$ is a dense $G_\delta$ subset of $\widetilde{X}$.
\end{proof}

Furstenberg introduced disjointness in topological dynamics in \cite{F67}.
We use the following version for Polish dynamical systems. 
For a Polish dynamical system $(X,T)$, a subset $A$ of $X$ is said to be \emph{$T$-invariant} if $TA\subset A$.
Let $(X,T)$ and $(Y,S)$ be two Polish dynamical systems.
A subset $G$ of $X\times Y$ is called a \emph{joining} of $(X,T)$ and $(Y,S)$ if it is closed and $T\times S$-invariant with $\overline{p_X(G)}=X$ and $\overline{p_Y(G)}=Y$, where $p_X$ and $p_Y$ denote the projection maps into $X$ and $Y$ respectively.
We say that $(X,T)$ is \emph{disjoint} from $(Y,S)$ if the only joining of $(X,T)$ and $(Y,S)$ is the whole space $X\times Y$.
Moreover, if both $(X,T)$ and $(Y,S)$ are compact dynamical systems, then the notion of disjointness coincides with that of Furstenberg in \cite{F67}. 

A compact dynamical system $(X,T)$ is said to be \emph{minimal} if the only nonempty closed $T$-invariant subset is the whole space $X$, or equivalently, $\Trans(T)=X$.
We next clarify the relationship between disjointness and transitive points.

\begin{proposition}\label{prop:compactification-disjoint}
    Let $(X,T)$ be a Polish dynamical system and let $(Y,S)$ be a minimal compact dynamical system.
    Consider the following assertions:
    \begin{enumerate}
        \item [(1)] $T$ is disjoint from $S$;
        \item [(2)] $\operatorname{Trans}(T)\times Y=\operatorname{Trans}(T\times S)$.
    \end{enumerate}
Then (1)$\Rightarrow$(2).
Moreover, if $T$ is transitive, then (1)$\Leftrightarrow$(2).
\end{proposition}
\begin{proof}
    (1)$\Rightarrow$(2). It suffices to show that $ \operatorname{Trans}(T)\times Y\subset \operatorname{Trans}(T\times S)$.
    Fix $x\in\operatorname{Trans}(T)$ and $y\in Y$.
    Then $G: =\overline{\operatorname{Orb}((x,y),T\times S)}$ is $T\times S$-invariant, and closed with $\overline{p_XG}=X$ and $\overline{p_YG}=Y$.
    So $G$ is a joining of $T$ and $S$.
    Since $T$ is disjoint from $S$, we have $G=X\times Y$ and then $(x,y)\in\operatorname{Trans}(T\times S)$.

    (2)$\Rightarrow$(1). Now assume that $T$ is transitive.
    Let $G$ be a joining of $T$ and $S$.
    Fix $x\in\operatorname{Trans}(T)$.
    Since $Y$ is compact, we have $X=\overline{p_{X}G}=p_{X}G$.
    Then there is $y\in Y$ with $(x,y)\in G$.
    By the hypothesis, we have $(x,y)\in\operatorname{Trans}(T\times S)$, which gives that $G=X\times Y$ since $G$ is closed and $T\times S$-invariant.
\end{proof}
Combining Propositions \ref{prop:transitive-compactification} and \ref{prop:compactification-disjoint}, we obtain that, for a transitive Polish dynamical system, disjointness from a given minimal compact dynamical system is equivalent to the disjointness of its compactifications from that minimal system.
\begin{proposition}\label{prop:disjoint-TS}
Let $(X,T)$ be a transitive Polish dynamical system with a compactification $(\widetilde{X}, \widetilde{T})$ and
    let $(Y,S)$ be a minimal compact dynamical system.
    Then $T$ is disjoint from $S$ if and only if $\widetilde{T}$ is disjoint from $S$.
\end{proposition}
\begin{proof}
    Suppose that $T$ is disjoint from $S$.
    Let $G$ be a joining of $\widetilde{T}$ and $S$.
    Since $Y$ is compact, we have $\operatorname{Trans}(T)\subset \widetilde{X}=\overline{p_{\widetilde{X}}G}=p_{\widetilde{X}}G$.
    Pick $x\in\operatorname{Trans}(T)$.
    Then there is $y\in Y$ with $(x,y)\in G$.
    By Proposition \ref{prop:compactification-disjoint}, $(x,y)\in\operatorname{Trans}(T\times S)$.
    Since $x\in X$, we have $\operatorname{Orb}((x,y), T\times S)=\operatorname{Orb}((x,y), \widetilde{T}\times S)$.
    So $\operatorname{Orb}((x,y), \widetilde{T}\times S)$ is dense in $\widetilde{X}\times Y$ and then $G=\widetilde{X}\times Y$, which gives that $\widetilde{T}$ is disjoint from $S$.

    Suppose that $\widetilde{T}$ is disjoint from $S$.
    By Propositions \ref{prop:transitive-compactification} and \ref{prop:compactification-disjoint},
    \begin{align*}
        \operatorname{Trans}(T)\times Y
        =&(\operatorname{Trans}(\widetilde{T})\cap X)\times Y
        =(\operatorname{Trans}(\widetilde{T})\times Y)\cap (X\times Y)\\
        =&\operatorname{Trans}(\widetilde{T}\times S)\cap (X\times Y)
        =\operatorname{Trans}(T\times S).
        \qedhere
    \end{align*}
    Since $T$ is transitive, by Proposition \ref{prop:transitive-compactification}, $\widetilde{T}$ is transitive. 
    By Proposition \ref{prop:compactification-disjoint}, $T$ is disjoint from $S$. 
\end{proof}

A disjointness fact arises from totally transitive compact dynamical systems with dense periodic points.
A Polish dynamical system $(X,T)$ is said to be \emph{totally transitive} if $T^p$ is transitive for every $p\in\mathbb N$.
A point $x\in X$ is called a \emph{periodic point} of $T$ if $T^px=x$ for some $p\in\mathbb N$.
We shall need the following disjointness result.

\begin{proposition}[{\cite[Theorem 2.2]{F67}}]\label{prop:compact-totally-transitive-disjoint-minimal}
    Let $(X,T)$ be a compact dynamical system.
    If $(X,T)$ is totally transitive with dense periodic points, then it is disjoint from every minimal compact dynamical system.
\end{proposition}

\begin{remark}\label{rem:totally-transitive-minimal}
Let $(X,T)$ be a Polish dynamical system with a compactification $(\widetilde{X}, \widetilde{T})$.
If $(X,T)$ is totally transitive with dense periodic points, then so is $(\widetilde{X}, \widetilde{T})$.
Therefore, by Propositions \ref{prop:compact-totally-transitive-disjoint-minimal} and \ref{prop:disjoint-TS}, $(X,T)$ is disjoint from every minimal compact dynamical system.
\end{remark}

Let $(X,T)$ and $(Y,S)$ be two Polish dynamical systems.
We say that $(X,T)$ is \emph{weakly disjoint} from $(Y,S)$ if the product system $(X\times Y, T\times S)$ is transitive.
By Proposition \ref{prop:transitive-compactification}, it follows that two Polish dynamical systems are weakly disjoint if and only if their compactifications are weakly disjoint.
Weak disjointness for hypercyclic operators is studied in \cite{LLR26}.

A compact dynamical system $(X,T)$ is said to be \emph{equicontinuous} if for every $\varepsilon>0$, there is $\delta>0$ such that for every pair $x,y\in X$ with $d(x,y)<\delta$, we have \[d(T^nx,T^ny)<\varepsilon,\quad n\in\mathbb N_0, \]
where $d$ is a compatible metric on $X$.
Notice that, since $X$ is compact, the above notion of equicontinuity is independent of the choice of a compatible metric on $X$.
Following \cite{HY12}, we say that a Polish dynamical system is \emph{weakly scattering} if it is weakly disjoint from every minimal equicontinuous compact dynamical system.
See \cite{HY12} for more details of weakly scattering compact dynamical systems.
Let $G$ be a compact metrizable group and let $g\in G$ be a topological generator.
Then the rotation $R_g\colon G\to G$, defined by $R_g(h)=gh$, is a minimal equicontinuous compact dynamical system.
Conversely, every minimal equicontinuous compact dynamical system is topologically conjugate to a rotation $(G,R_g)$ on a compact metrizable group $G$ by a topological generator $g\in G$, see e.g.\ \cite{F63}.
Therefore, a Polish dynamical system is weakly scattering if and only if it is weakly disjoint from every such rotation system $(G,R_g)$.

Let $(X,T)$ be a compact dynamical system and let $d$ be a compatible metric on $X$.
We say that $(X,T)$ is \emph{distal} if for every pair $x,y\in X$ with $x\neq y$, we have
\[
    \inf_{n\in\mathbb N_0} d(T^nx,T^ny)>0.
\]
Notice that, since $X$ is compact, the above notion of distality is independent of the choice of a compatible metric on $X$.
The above rotation $(G,R_g)$ is an example of a distal system.
Distal compact dynamical systems form a classical class in topological dynamics, see \cite{F63}.
In a distal compact dynamical system the transformation is a homeomorphism, and the orbit closure of each point is minimal \cite{F63}.
In particular, a distal system with a transitive point is minimal.
Subsystems (i.e.\ closed invariant subsets) and finite or countable products of distal systems are distal. 
Let $(X,T)$ and $(Y,S)$ be two compact dynamical systems.
If there exists a continuous surjection $\Phi\colon X\to Y$ with $\Phi T=S\Phi$, then we say that $(Y,S)$ is a \emph{factor} of $(X,T)$.
Every factor of a distal compact dynamical system is also distal, see \cite[Theorem 3.3]{F63}.
In \cite{HY12}, Huang and Ye proved the following disjointness fact concerning weakly scattering systems and minimal distal systems.

\begin{proposition}[{\cite[Theorem 2.14]{HY12}}]\label{prop:HY12}
Let $(X,T)$ be a compact dynamical system.
If $(X,T)$ is weakly scattering, then it is disjoint from every minimal distal compact dynamical system.
\end{proposition}

\begin{remark}\label{rem:weak-scattering-disjoint-distal}
Let $(X,T)$ be a Polish dynamical system with a compactification $(\widetilde{X}, \widetilde{T})$.
    If $(X,T)$ is weakly scattering, then so is $(\widetilde{X}, \widetilde{T})$.
    By Propositions \ref{prop:HY12} and \ref{prop:disjoint-TS}, $(X,T)$ is disjoint from every minimal distal compact dynamical system.
\end{remark}

\subsection{Hypercyclic operators}
Let $T$ be an operator on a complex F-space $X$.
A vector $x\in X$ is said to be \emph{hypercyclic} for $T$ if $N_T(x,U): =\{n\in\mathbb N_0\colon T^nx\in U\}$ is nonempty for every nonempty open subset $U$ of $X$, or equivalently, its orbit $\Orb(x,T):=\{T^nx\colon n\in\mathbb N_0\}$ is dense in $X$. 
The set of hypercyclic vectors is denoted by $\HC(T)$. 
We say that $T$ is \emph{hypercyclic} if $\HC(T)\neq\emptyset$. 
In this case, $X$ is separable, so $(X,T)$ is a Polish dynamical system. 
By Theorem \ref{thm:birkhoff}, $(X,T)$ is transitive and $\HC(T)=\Trans(T)$ is a dense $G_\delta$ subset of $X$. 
Hypercyclicity with respect to Furstenberg families has been systematically studied in linear dynamics, see, for instance, \cite{BGE07}, \cite{BGE18} and \cite{BGLP22}.
Let $\cf$ be a family.
We denote 
\[\HC_{\cf}(T):=\left\{x\in X\colon N_T(x,U)\in\cf, \forall\text{ nonempty open }U\subset X\right\}\]
and call the elements of $\HC_{\cf}(T)$ the \emph{$\cf$-hypercyclic} vectors of $T$.
We say that $T$ is \emph{$\cf$-hypercyclic} if $\HCF(T)\neq\emptyset$.

Let $\bs\lambda=(\lambda_n)_n$ be a complex sequence. 
For $x\in X$ and $U\subset X$, write 
\[
    N_T^{\bs\lambda}(x,U)=\{n\in\mathbb N_0\colon \lambda_nT^n x\in U\}
\] 
and denote 
\[\HCF^{\bs\lambda}(T): =\left\{x\in X\colon N_T^{\bs\lambda}(x,U)\in\cf,\  \forall \text{ nonempty open subset }U\text{ of }X\right\}. 
\]
We simply write $\HC^{\bs\lambda}(T)=\HC^{\bs\lambda}_{\cf_{\mathrm{inf}}}(T)$.
One observes that \[\HC^{\bs\lambda}(T)=\left\{x\in X\colon \{\lambda_nT^nx\colon n\in\mathbb N_0\} \text{ is dense in }X\right\}. \]
It can be seen that $\HCF^{\bs\lambda}(T)=\HCF^{c\bs\lambda}(T)$ for every $c\in\mathbb C\setminus\{0\}$. 
The next observation allows finite changes and small relative perturbations. 

\begin{proposition}\label{prop:sequence-finite-mod}
Let $T$ be an operator on a complex F-space $X$, let $\cf$ be a family and let $\bs\lambda$ and $\bs\mu$ be two complex sequences.
If $\mu_n/\lambda_n\to 1$ as $n\to\infty$ along the indices $n$ for which $\lambda_n\neq 0$, then $\HCF^{\bs\lambda}(T)\subset \HCF^{\bs\mu}(T)$.
\end{proposition}
\begin{proof}
We may assume that $X$ is not a singleton.
Fix $x\in\HCF^{\bs\lambda}(T)$ and a nonempty open set $U\subset X$.
Choose a nonempty open set $U_0\subset U$ with $0\notin U_0$.
Continuity of scalar multiplication gives a neighborhood $D$ of $1$ and a nonempty open set $V\subset U_0$ such that $D\cdot V\subset U_0$.
Choose $N$ such that $\mu_n/\lambda_n\in D$ whenever $n>N$ and $\lambda_n\ne0$.
The open set
\[
W=V\setminus\{\lambda_lT^lx\colon 0\leq l\leq N\}
\]
is nonempty.
For $n\in N_T^{\bs\lambda}(x,W)$, we have $n>N$ and $\lambda_n\ne0$.
Therefore
\[
\mu_nT^nx=(\mu_n/\lambda_n)\lambda_nT^nx\in D\cdot V\subset U.
\]
Thus $N_T^{\bs\lambda}(x,W)\subset N_T^{\bs\mu}(x,U)$.
The left-hand side belongs to $\cf$, so the right-hand side does as well.
\end{proof}

\begin{remark}\label{rem:sequence-finite-mod}
Proposition \ref{prop:sequence-finite-mod} has the following consequences.
\begin{enumerate}
\item Replacing zero terms of $\bs\lambda$ by arbitrary scalars can only enlarge $\HCF^{\bs\lambda}(T)$ since the relative ratio on the nonzero terms is exactly $1$.
\item Changing finitely many terms of the weight $\bs\lambda$ leaves $\HCF^{\bs\lambda}(T)$ unchanged, by applying the proposition in both directions.
\item If $\mu_n\to1$ as $n\to\infty$, then $\HCF^{\bs\lambda}(T)=\HCF^{\bs\lambda\bs\mu}(T)$, where $\bs\lambda\bs\mu=(\lambda_n\mu_n)_n$.
Again both relative ratios tend to $1$ on the relevant nonzero terms.
\end{enumerate}
\end{remark}

Let $c_0$ be the space of complex sequences tending to zero. 
Weighted hypercyclicity is stable under $c_0$-perturbations when the weight is bounded away from zero. 
More precisely, we have the following. 

\begin{proposition}
Let $T$ be an operator on a complex F-space $X$, let $\cf$ be a family, let $\bs\lambda$ be a complex sequence and let $\bs\mu\in c_0$. 
If the nonzero values of $\bs\lambda$ are bounded away from zero, then $\HCF^{\bs\lambda}(T)\subset \HCF^{\bs\lambda+\bs\mu}(T)$.
Moreover, if $\{n\in\mathbb N_0\colon \lambda_n=0\}$ is finite, then $\HCF^{\bs\lambda}(T)=\HCF^{\bs\lambda+\bs\mu}(T)$.
\end{proposition}
\begin{proof}
Choose $a>0$ with $|\lambda_n|\geq a$ whenever $\lambda_n\ne0$.
Then
\[
\left|\frac{\lambda_n+\mu_n}{\lambda_n}-1\right|
\leq a^{-1}|\mu_n|\to0
\quad(\lambda_n\ne0,n\to\infty).
\]
The first conclusion follows from Proposition \ref{prop:sequence-finite-mod}. 

If $\{n\in\mathbb N_0\colon \lambda_n=0\}$ is finite, then for all sufficiently large $n$,
$|\lambda_n|\geq a$ and $|\lambda_n+\mu_n|\geq a/2$.
Consequently,
\[
\left|\frac{\lambda_n}{\lambda_n+\mu_n}-1\right|
\leq 2a^{-1}|\mu_n|\to0,n\to\infty.
\]
A second application of Proposition \ref{prop:sequence-finite-mod} yields
$\HCF^{\bs\lambda+\bs\mu}(T)\subset\HCF^{\bs\lambda}(T)$.
Combining the two inclusions, we then obtain the equality.
\end{proof}

\section{From weighted hypercyclicity to ordinary hypercyclicity}
In this section, we study the removal of scalar weights. 
We first prove Theorem \ref{thm:removal}, and then establish a finer characterization for fixed families satisfying additional combinatorial conditions. 
The necessity of both bounds in Theorem \ref{thm:removal} follows by coding a dense set into disjoint blocks for the backward shift $B$ on $\mathbb C^{\mathbb N_0}$.

\begin{proposition}\label{prop:HWS-bd}
Let $\bs\lambda$ be a complex sequence.
If $\HC^{\bs\lambda}(B)\subset \HC(B)$, then the nonzero values of $\bs\lambda$ are bounded and bounded away from zero, where $B$ is the backward shift on $\mathbb C^{\mathbb N_0}$, that is, $B(z_0,z_1,\dotsc)=(z_1,z_2,\dotsc)$. 
\end{proposition}

\begin{proof}
Let $(y^{(k)})_k$ be a sequence of nonzero finite-support vectors which is dense in $\mathbb C^{\mathbb N_0}$.
For each $k\in\mathbb N_0$, choose
\[
M_k>\max\Biggl(
2\Bigl(\max_j|y_j^{(k)}|+1\Bigr),
\frac{3}{\min\{|y_j^{(k)}|\colon y_j^{(k)}\neq0\}}
\Biggr).
\]
Denote $l_k=\max\{t\in\mathbb N_0:y_t^{(k)}\neq0\}$ for each $k\in\mathbb N_0$.
Assume on the contrary that the nonzero values of $\bs\lambda$ are not bounded or not bounded away from zero.
Then $E_M:=\{n\in\mathbb N_0:|\lambda_n|>M\text{ or }0<|\lambda_n|<M^{-1}\}$ is infinite for every $M\geq1$.
Hence we can choose a strictly increasing sequence $(p_k)_k$ such that $p_k\in E_{M_k}$ and $p_{k+1}>p_k+l_k$ for every $k\in\mathbb N_0$.
Then the intervals $([p_k,p_k+l_k])_k$ are pairwise disjoint.
Define $x\in\mathbb C^{\mathbb N_0}$ by
\[
x_n=
\begin{cases}
\dfrac{y_j^{(k)}}{\lambda_{p_k}},
& n=p_k+j\text{ for some }k\in\mathbb N_0
\text{ and }j\leq l_k,\\[2mm]
0,&\text{else}.
\end{cases}
\]

Let $U$ be a nonempty open subset of $\mathbb C^{\mathbb N_0}$.
Then there exist $u\in U$, $s\in\mathbb N_0$ and $\varepsilon>0$ such that $V:=\{v\in\mathbb C^{\mathbb N_0}:|v_j-u_j|<\varepsilon\text{ for }j=0,\ldots,s\}\subset U$.
Let $W=V\cap\{v\in\mathbb C^{\mathbb N_0}:v_{s+1}\neq0\}$, which is nonempty and open.
Then there exists $k\in\mathbb N_0$ with $y^{(k)}\in W$, which gives that $l_k\geq s+1$.
For every $j\leq s$, we have
\[
|(\lambda_{p_k}B^{p_k}x)_j-u_j|
=|\lambda_{p_k}x_{p_k+j}-u_j|
=|y_j^{(k)}-u_j|<\varepsilon.
\]
This gives that $\lambda_{p_k}B^{p_k}x\in V\subset U$.
Hence $x\in\HC^{\bs\lambda}(B)$.

On the other hand, for every $n\in\mathbb N_0$ with $x_n\neq0$, there exist $k\in\mathbb N_0$ and $j\leq l_k$ such that $x_n=y_j^{(k)}/\lambda_{p_k}$.
Since $p_k\in E_{M_k}$, we have either $|\lambda_{p_k}|>M_k$ or $0<|\lambda_{p_k}|<M_k^{-1}$.
In the first case,
\[
|x_n|
=\biggl|\frac{y_j^{(k)}}{\lambda_{p_k}}\biggr|
<\frac{\max_j\bigl|y_j^{(k)}\bigr|}{M_k}
<\frac12.
\]
In the second case,
\[
|x_n|
=\biggl|\frac{y_j^{(k)}}{\lambda_{p_k}}\biggr|
>M_k\min\bigl\{|y_j^{(k)}|\colon y_j^{(k)}\neq0\bigr\}
>3.
\]
Thus, for every $n\in\mathbb N_0$, we have $|x_n|<1/2$ or $|x_n|>3$.
Consequently, $x\notin\HC(B)$, which contradicts the hypothesis.
\end{proof}

Let $T$ be an operator on a complex F-space $X$ and let $\Gamma\subset\mathbb C$.
For $x\in X$ and $U\subset X$, define
\[
N_T(\Gamma x,U)=\{n\in\mathbb{N}_0\colon \gamma T^nx\in U\text{ for some }\gamma\in\Gamma\}.
\]
A vector $x$ is \emph{$\Gamma$-supercyclic} if $N_T(\Gamma x,U)$ is nonempty for every nonempty open subset $U$ of $X$. 
We denote the set of these vectors by $\operatorname{SC}^{\Gamma}(T)$.
The term \emph{supercyclic} refers to $\Gamma=\mathbb C$.

The scalar-set theorem of Charpentier, Ernst and Menet is stated for Banach spaces in \cite[Theorem A]{CEM16}. 
Its Fr\'echet-space version is noted in \cite[Remark 1.2]{CEM16}.
We shall use the following extension to complex F-spaces.

\begin{proposition}\label{thm:CEM}
Let $\Gamma$ be a nonempty subset of $\mathbb{C}$. 
Then the following assertions are equivalent. 
\begin{enumerate}
    \item $\Gamma\setminus\{0\}$ is nonempty, bounded, and bounded away from zero;
    \item for every infinite-dimensional complex Banach space $X$ and every operator $T\colon X\to X$, $\operatorname{SC}^{\Gamma}(T)=\HC(T)$;
    \item for every complex F-space $X$ and every operator $T\colon X\to X$, $\operatorname{SC}^{\Gamma}(T)=\HC(T)$. 
\end{enumerate}
\end{proposition}
\begin{proof}
The equivalence of (1) and (2) is exactly \cite[Theorem A]{CEM16} and the implication (3)$\Rightarrow$(2) is clear. 
It remains to show that (1)$\Rightarrow$(3). 
Let $T$ be an operator on a complex F-space $X$ which is not a singleton. 
Since $\Gamma\setminus\{0\}$ is nonempty, we have $\HC(T)\subset\operatorname{SC}^{\Gamma}(T)$. 
In the following, we assume that $\operatorname{SC}^{\Gamma}(T)\neq \emptyset$, fix $x\in \operatorname{SC}^{\Gamma}(T)$ and shall show that $x$ is hypercyclic for $T$. 
Assume that $\dim X<\infty$. 
Since $T$ is supercyclic on a complex F-space $X$, we obtain that $\dim X=1$. 
Hence $T=\lambda I$ for some $\lambda\in\mathbb C$. 
Let $a,b>0$ be such that $\Gamma\setminus\{0\}\subset\{z\in\mathbb C\colon a\leq |z|\leq b\}$. 
For every $n\in\mathbb N_0$ and $\gamma\in\Gamma\setminus\{0\}$, if $|\lambda|\leq 1$, we have 
$|\gamma\lambda^nx|\leq b|\lambda|^n|x|\leq b|x|$. 
If $|\lambda|> 1$, we have 
$|\gamma\lambda^nx|\geq a|x|$. 
That is, $\{\gamma T^nx\colon \gamma\in\Gamma,n\in\mathbb N_0\}$ is not dense in $X$, a contradiction. 
So we have $\dim X=\infty$. 
We need the following technical claim. 

\medskip 
\noindent \textbf{Claim}: For every nonzero complex polynomial $P$, $P(T)$ has dense range, and $P(T)x$ is $\Gamma$-supercyclic. 
\begin{proof}[Proof of the Claim]
    Let $M=\overline{P(T)X}$, which is a closed subspace of $X$. 
Assume that $M\neq X$. 
Consider the nonzero-dimensional quotient space $Q: =X/M$. 
Let $q\colon X\to Q$ be the quotient map. 
Since $M$ is $T$-invariant, $T$ induces a map $\widehat{T}\colon Q\to Q$, where 
\[
\widehat{T}(q(y))=q(Ty),\quad y\in X. 
\]
Since $\operatorname{Orb}(\Gamma x,T): =\bigl\{\gamma T^nx\colon \gamma\in \Gamma, n\in\mathbb N_0\bigr\}$ is dense in $X$ and $q$ is a continuous surjection, we have that $D: =q(\operatorname{Orb}(\Gamma x,T))=\bigl\{\gamma\widehat{T}^nq(x)\colon \gamma\in \Gamma, n\in\mathbb N_0\bigr\}$ is dense in $Q$. 
For each $y\in X$, since $P(T)y\in P(T)X\subset M$, we have $P(\widehat{T})(q(y))=q(P(T)y)=0\in Q$. 
Without loss of generality, write $P(z)=a_dz^d+\dotsb +a_0$ with $a_d\neq 0$ and $d\geq 1$. 
Since $P(\widehat{T})=0$\ (the zero map), we have $P(\widehat{T})=a_d\widehat{T}^d+\dotsb +a_0I=0$. 
Since $a_d\neq 0$, we have $\widehat{T}^d\in\operatorname{span}\{\widehat{T}^i\colon 0\leq i\leq d-1\}$. 
In particular, for each $n\in\mathbb N_0$, 
\[
    \widehat{T}^nq(x)\in E: =\operatorname{span}\{\widehat{T}^i(q(x))\colon 0\leq i\leq d-1\}. 
\]
So the dense subset $D$ of $Q$ is completely contained in $E$. 
Since $E$ is of finite dimension, it is closed and then $Q=E$. 
Now we consider the adjoint $\widehat{T}^*\colon Q^*\to Q^*$. 
Since $\widehat{T}^*$ admits an eigenvector, there is a nonzero $f\in Q^*$ and $\alpha\in\mathbb C$ with $\widehat{T}^*f=\alpha f$. 
In particular, for each $n\in\mathbb N_0$ and $\gamma\in \Gamma$, 
\[
    f(\gamma\widehat{T}^nq(x))=\gamma \alpha^nf(q(x)). \tag{\raisebox{-0.47ex}{*}}\label{star}
\]
If $f(q(x))=0$, since $f$ is continuous, we have $f|_Q=0$. 
That is, $f$ is the zero functional, which is impossible. 
So $f(q(x))\neq 0$. 
Since $f$ is nonzero, the left-hand side of \eqref{star} equals $f(D)$, which is dense in $\mathbb C$. 
So the right-hand side of \eqref{star} is also dense in $\mathbb C$. 
Since $\Gamma\setminus\{0\}$ is bounded and bounded away from zero, there exist $r_1,r_2>0$ such that $r_1\leq |\gamma|\leq r_2$ for every $\gamma\in\Gamma\setminus\{0\}$. 

If $|\alpha|\leq 1$, then for every $n\in\mathbb N_0$ and $\gamma\in \Gamma\setminus\{0\}$, 
\[|\gamma \alpha^nf(q(x))|=|\gamma|\cdot |\alpha^n|\cdot |f(q(x))|\leq r_2 |f(q(x))|. \]
If $|\alpha|>1$, then for every 
$n\in\mathbb N_0$ and $\gamma\in \Gamma\setminus\{0\}$, 
\[|\gamma \alpha^nf(q(x))|=|\gamma|\cdot |\alpha^n|\cdot |f(q(x))|\geq r_1 |f(q(x))|. \]
So the right-hand side of (*) can never be dense, which is a contradiction. 
So $M=X$. 
It is then easy to see that $P(T)x$ is also $\Gamma$-supercyclic. 
\end{proof}
We replace Lemmas 3.5 and 3.6 of \cite{CEM16} with the above claim. 
All the other ingredients used in the proof of Theorem 3.2 of \cite{CEM16} remain valid for complex F-spaces, 
with the norm estimates replaced, when necessary, by their counterparts for a compatible metric. 
Therefore, the same argument applies and yields that $x$ is a hypercyclic vector of $T$. 
\end{proof}

For a family $\cf$, a vector $x$ is \emph{$(\cf,\Gamma)$-supercyclic} if $N_T(\Gamma x,U)\in\cf$ for every nonempty open set $U\subset X$.
Write $\operatorname{SC}^{\Gamma}_{\cf}(T)$ for the set of such vectors.
This notion was studied by Alves, Botelho and F\'avaro \cite{ABF26}.
Proposition 5.9 of \cite{ABF26} gives an annular scalar-removal result on Banach spaces for right-invariant hypercyclicity sets whose complements are closed under finite unions.
Here the term \emph{hypercyclicity set} requires a sequence of pairwise disjoint members $A_k\in\cf$ such that distinct $n\in A_k$ and $m\in A_l$ satisfy $|n-m|\geq\max\{k,l\}$.
The following proposition uses the CuSP property instead and applies on arbitrary complex F-spaces.
In particular, it covers $\cf_{\mathrm{pld}}$, whose complement is not closed under finite unions.

\begin{proposition}\label{prop:annular-removal}
Let $\Gamma$ be a subset of $\C$.
If $\Gamma\setminus\{0\}$ is nonempty, bounded, and bounded away from zero, then for every complex F-space $X$, every operator $T\colon X\to X$ and every family $\cf$ with the CuSP property, we have
\[\operatorname{SC}^{\Gamma}_{\cf}(T)=\HCF(T). \]
\end{proposition}

\begin{proof}
We may assume that $X$ is not a singleton.
Fix $\gamma_0\in\Gamma\setminus\{0\}$.
For $x\in\HCF(T)$ and nonempty open $U\subset X$,
\[
N_T(x,\gamma_0^{-1}U)\subset N_T(\Gamma x,U).
\]
This proves $\HCF(T)\subset\operatorname{SC}^{\Gamma}_{\cf}(T)$.

Conversely, fix $x\in\operatorname{SC}^{\Gamma}_{\cf}(T)$ and nonempty open $U\subset X$.
The scalar orbit is dense, so Proposition \ref{thm:CEM} gives $x\in\HC(T)$. 
In particular $x\ne0$.
Choose a nonempty open set $V\subset U$ and $\varepsilon>0$ such that $D(1,\varepsilon)V\subset U$, where $D(a,r)=\{z\in\mathbb C\colon |z-a|<r\}$.
Let $c>0$ be a lower bound for $|\gamma|$ on $\Gamma\setminus\{0\}$, and put $K=\overline{\Gamma\setminus\{0\}}$, which is compact. 
Choose $\mu_1,\dotsc,\mu_q\in K$ with $K\subset\bigcup_jD(\mu_j,c\varepsilon)$.
Since $x$ is hypercyclic, there exist integers $r_j\in\mathbb{N}_0$ with $T^{r_j}x\in\mu_jV$, $j=1,2,\dotsc,q$. 
Choose a neighborhood $W$ of $x$ with $0\notin W$ such that $T^{r_j}W\subset\mu_jV$ for $j=1,2,\dotsc,q$.

Put $A=N_T(\Gamma x,W)\in\cf$.
For each $n\in A$, there exists $\gamma_n\in\Gamma\setminus\{0\}$ with $\gamma_nT^nx\in W$. 
For $j=1,2,\dotsc,q$, let $A_j=\{n\in A\colon \gamma_n\in D(\mu_j,c\varepsilon)\}$.
Then $A=\bigcup_{j=1}^qA_j$.
Fix $j\in\{1,2,\dotsc,q\}$. 
For $n\in A_j$ we have $|\gamma_n^{-1}\mu_j-1|<\varepsilon$, and hence
\[
T^{n+r_j}x
=\gamma_n^{-1}T^{r_j}(\gamma_nT^nx)
\in\gamma_n^{-1}\mu_jV\subset U.
\]
Thus $\bigcup_{j=1}^q(r_j+A_j)\subset N_T(x,U)$.
By the CuSP property, we obtain $N_T(x,U)\in\cf$. 
Hence $x\in\HCF(T)$. 
\end{proof}

We now complete the proof of Theorem \ref{thm:removal}. 
\begin{proof}[Proof of Theorem \ref{thm:removal}]
(1)$\Rightarrow$(2) follows from Proposition \ref{prop:HWS-bd}. 

(2)$\Rightarrow$(1). If $\bs\lambda=\bs 0$, then the result is trivial. 
Assume that $\bs\lambda\neq\bs 0$ and put $\Gamma=\{\lambda_n\colon n\in\mathbb{N}_0\}$.
It is observed that $\HC^{\bs\lambda}(T)\subset \operatorname{SC}^{\Gamma}(T)$. 
Since $\Gamma\setminus\{0\}$ is nonempty, bounded and bounded away from zero, by Proposition \ref{thm:CEM}, $\operatorname{SC}^{\Gamma}(T)=\HC(T)$. 
Hence $\HC^{\bs\lambda}(T)\subset\HC(T)$. 

Moreover, suppose that the nonzero values of $\bs\lambda$ are bounded and bounded away from zero, and let $\cf$ be a family with the CuSP property.
Again we assume that $\bs\lambda\neq \bs 0$ and denote $\Gamma=\{\lambda_n\colon n\in\mathbb{N}_0\}$.
It is observed that $\HCF^{\bs\lambda}(T)\subset \operatorname{SC}_{\cf}^{\Gamma}(T)$. 
Since $\Gamma\setminus\{0\}$ is nonempty, bounded and bounded away from zero, by Proposition \ref{prop:annular-removal}, $\operatorname{SC}_{\cf}^{\Gamma}(T)=\HCF(T)$. 
Hence $\HCF^{\bs\lambda}(T)\subset\HCF(T)$. 
\end{proof}

Theorem \ref{thm:removal} also has a version for individual vectors when the exceptional weighted orbit is nowhere dense.
A convergent exceptional orbit is one example.

\begin{remark}
    Let $\bs\lambda$ be a complex sequence, let $T$ be an operator on a complex F-space $X$ which is not a singleton and let $\cf$ be a family with the CuSP property.
    For $x\in \HCF^{\bs\lambda}(T)$, if there exists $E\subset \mathbb N_0$ such that $\{\lambda_n\colon n\in E^c\}\setminus\{0\}$ is nonempty, bounded and bounded away from zero, and $K_x=\overline{\{\lambda_nT^nx\colon n\in E\}}$ has empty interior, then for every nonempty open subset $U$ of $X$, the set $U\setminus K_x$ is also nonempty and open.
    Since $N_T^{\bs\lambda}(x,U\setminus K_x)\subset N_T(\{\lambda_n\colon n\in E^c\}x,U)$, we have that $x\in\operatorname{SC}_{\cf}^{\{\lambda_n\colon n\in E^c\}}(T)$.
    By Proposition \ref{prop:annular-removal}, $x\in\HCF(T)$.
\end{remark}

For the family $\cf_{\mathrm{inf}}$, Theorem \ref{thm:removal} imposes uniform bounds on all nonzero weights. 
For a fixed family, these bounds can be relaxed to a condition involving its dual family. 
The following theorem gives a characterization under a suitable separation
assumption on the family. 

\begin{theorem}\label{thm:F-removal}
    Assume that a family $\cf$ has the CuSP property and Ramsey property, and for every sequence $(E_k)_{k\geq 0}$ in $\cf$ and every sequence $(q_k)_{k\geq 0}$ of $\mathbb N$, there exist a sequence $(F_k)_{k\geq 0}$ in $\cf$ such that 
    \begin{enumerate}
        \item[(a)] $F_k\subset E_k$ for every $k\in\mathbb N_0$ and $F_{k_1}\cap F_{k_2}=\emptyset$ for $k_1\neq k_2$; 
        \item[(b)] $|i-j|\geq\max(q_{k_1},q_{k_2})$ for every $k_1,k_2\in\mathbb N_0$ and $i\in F_{k_1}, j\in F_{k_2}$ with $i\neq j$. 
    \end{enumerate}
    Then for every complex sequence  $\bs\lambda$, the following assertions are equivalent. 
    \begin{enumerate}
        \item $\HCF^{\bs\lambda}(T)\subset \HCF(T)$ for every complex F-space $X$ and every operator $T\colon X\to X$; 
        \item there exist $0<a\leq b$ such that $\{n\in\mathbb N_0\colon \lambda_n=0\text{ or }a\leq |\lambda_n|\leq b\}\in\cf^*$. 
    \end{enumerate}
\end{theorem}

\begin{proof}
    (2)$\Rightarrow$(1). 
    Put $\Gamma=\{z\in\mathbb C\colon a\leq |z|\leq b\}\cup\{0\}$. 
    Let $x\in \HCF^{\bs\lambda}(T)$ and let $U$ be a nonempty open subset of $X$. 
    Denote $A=\{n\in\mathbb N_0\colon \lambda_n=0\text{ or }a\leq |\lambda_n|\leq b\}$. 
    For every $n\in N_T^{\bs\lambda}(x,U)\cap A$, we have $\lambda_nT^nx\in U$ and $\lambda_n\in\Gamma$. 
    This gives that $N_T^{\bs\lambda}(x,U)\cap A\subset N_T(\Gamma x,U)$. 
    Since $\cf$ has the Ramsey property, we have $N_T^{\bs\lambda}(x,U)\cap A\in\cf$. 
    Hence $x\in\operatorname{SC}^{\Gamma}_{\cf}(T)$. 
    By Proposition \ref{prop:annular-removal}, since $\cf$ has the CuSP property and $\Gamma\setminus\{0\}$ is bounded and bounded away from zero, $x\in\HCF(T)$. 

    (1)$\Rightarrow$(2). Let $(y^{(k)})_k$ be a sequence of nonzero finite-support vectors which is dense in $\mathbb C^{\mathbb N_0} $. 
    For each $k\in\mathbb N_0$, choose 
\[
M_k>\max\Biggl(
2\Bigl(\max_j|y_j^{(k)}|+1\Bigr),
\frac{3}{\min\bigl\{|y_j^{(k)}|\colon y_j^{(k)}\neq0\bigr\}}
\Biggr).
\]
    Denote $l_k=\max\{t\in\mathbb N_0\colon y_t^{(k)}\neq 0\}$ for each $k\in\mathbb N_0$. 
    Assume on the contrary that (2) fails. 
    Then for every $M\geq 1$, we have $A_M=\{n\in\mathbb N_0\colon \lambda_n=0\text{ or }M^{-1}\leq |\lambda_n|\leq M\}\notin\cf^*$. 
    Then $E_M: =\mathbb N_0\setminus A_M\in\cf$ for every $M\geq 1$. 
    Considering $(E_{M_k})_k$ and $(l_k+1)_k$, by the hypothesis, we obtain a sequence $(F_k)_k$ in $\cf$ which is pairwise disjoint such that $F_k\subset E_{M_k}$, and $|i-j|\geq\max(l_{k_1}+1,l_{k_2}+1)$ for every $k_1,k_2\in\mathbb N_0$ and $i\in F_{k_1}, j\in F_{k_2}$ with $i\neq j$. 
    Then the intervals $([m,m+l_k])_{m\in F_k,k\in\mathbb N_0}$ are pairwise disjoint. 
    Define $x\in\mathbb C^{\mathbb N_0}$ by 
    \[
        x_n=\begin{cases}
            \frac{y_j^{(k)}}{\lambda_m}, &n=m+j\text{ for some }k\in\mathbb N_0,m\in F_k\text{ and }j\leq l_k, \\
            0, &\text{else}. 
        \end{cases}
    \]

Let $U$ be a nonempty open subset of $\mathbb C^{\mathbb N_0} $. 
Then there exist $u\in U$, $s\in\mathbb N_0$ and $\varepsilon>0$ such that $V\colon =\{v\in\mathbb C^{\mathbb N_0} \colon |v_j-u_j|<\varepsilon\text{ for }j=0,1,\dotsc,s\}\subset U$. 
Let $W=V\cap \{v\in\mathbb C^{\mathbb N_0} \colon v_{s+1}\neq 0\}$, which is nonempty and open. 
Then there exists $k\in\mathbb N_0$ with $y^{(k)}\in W$, which gives that $l_k\geq s+1$. 
    For $m\in F_k$, then we have \[|(\lambda_{m}B^{m}x)_j-u_j|=|\lambda_{m}x_{m+j}-u_j|=|y^{(k)}_j-u_j|<\varepsilon, \quad j\leq s. \]
    This gives that $F_k\subset N_B^{\bs\lambda}(x,V)\subset N_B^{\bs\lambda}(x,U)$. 
    Hence $x\in\HCF^{\bs\lambda}(B)$. 

    On the other hand, for every $n\in\mathbb N_0$ with $x_n\neq 0$,
there exist $k\in\mathbb N_0$, $m\in F_k$, and $j\le l_k$ such that $x_n=y_j^{(k)}/\lambda_m$. 
Since $m\in F_k\subset E_{M_k}$, we have either
$|\lambda_m|>M_k$ or $0<|\lambda_m|<M_k^{-1}$.
In the first case, we have
\[
|x_n|
=\biggl|\frac{y_j^{(k)}}{\lambda_m}\biggr|
<\frac{\max_j|y_j^{(k)}|}{M_k}
<\frac12.
\]
In the second case, we have
\[
|x_n|
=\biggl|\frac{y_j^{(k)}}{\lambda_m}\biggr|
>M_k\min\bigl\{|y_j^{(k)}|\colon y_j^{(k)}\neq 0\bigr\}
>3.
\]
Thus, for every $n\in\mathbb N_0$, we have
\[
|x_n|<\frac12
\quad\text{or}\quad
|x_n|>3.
\]
Hence $x\notin\HC(B)$, and thus
$x\notin\mathrm{HC}_{\mathcal F}(B)$,
which contradicts (1). 
Therefore (2) holds.
\end{proof}

\begin{corollary}\label{coro:F-removal}
    Let $\bs\lambda$ be a complex sequence. 
    Then for $\cf=\cf_{\mathrm{inf}},\cf_{\mathrm{pubd}},\cf_{\mathrm{pud}},\cf_{\mathrm{ps}}$, the following assertions are equivalent. 
 \begin{enumerate}
     \item $\HCF^{\bs\lambda}(T)\subset \HCF(T)$ for every complex F-space $X$ and every operator $T\colon X\to X$; 
     \item there exist $0<a\leq b$ such that $\{n\in\mathbb N_0\colon \lambda_n=0\text{ or }a\leq |\lambda_n|\leq b\}\in \cf^*$. 
 \end{enumerate}
 \end{corollary}

\begin{proof}
The families $\cf_{\mathrm{inf}}$, $\cf_{\mathrm{pubd}}$, $\cf_{\mathrm{pud}}$, and $\cf_{\mathrm{ps}}$ have both the CuSP and Ramsey properties. 
Moreover, they satisfy the separation condition in Theorem \ref{thm:F-removal}, as verified below. 
The conclusion now follows from Theorem \ref{thm:F-removal}. 

Fix sequences $(E_k)_k$ in the family under consideration and $(q_k)_k$ in $\mathbb N$.
Choose $(k_n)_n$ in which every nonnegative integer occurs infinitely often, and put $Q_n=\max_{0\leq t\leq n}q_{k_t}$.

\textbf{Case 1: $\cf_{\mathrm{inf}}$.}
Recursively choose a strictly increasing sequence $(a_n)_n$ with $a_n\in E_{k_n}$ and $a_n-a_{n-1}\geq Q_n$ for $n\geq1$.
Then $F_k=\{a_n\colon k_n=k\}$ ($k\in\mathbb N_0$) are desired.

\textbf{Case 2: $\cf_{\mathrm{pubd}}$.}
Since $\cf_{\mathrm{pubd}}$ has the Ramsey property, for each $k$ there exists $r_k\in\{0,\ldots,q_k-1\}$ such that $P_k:=\{n\in E_k\colon n\equiv r_k\pmod{q_k}\}\in\cf_{\mathrm{pubd}}$.
Write $\delta_k=\overline{\operatorname{Bd}}(P_k)>0$.
Recursively choose finite integer intervals $I_n=[a_n,b_n]\cap\mathbb N_0$ such that
\[
|I_n|\geq n+1,\quad
|P_{k_n}\cap I_n|\geq\frac{\delta_{k_n}}2|I_n|,
\quad
a_n-b_{n-1}>Q_n\quad(n\geq1).
\]
Put $B_n=P_{k_n}\cap I_n$ and
\[
F_k=\bigcup_{\{n\colon k_n=k\}}B_n.
\]
For each $k$, infinitely many $n$ satisfy $k_n=k$, and along these indices $|I_n|\to\infty$.
Consequently,
\[
\overline{\operatorname{Bd}}(F_k)
\geq\limsup_{\substack{n\to\infty\\k_n=k}}
\frac{|F_k\cap I_n|}{|I_n|}
\geq\frac{\delta_k}{2}>0.
\]
The $I_n$ are pairwise disjoint, hence so are the $F_k$.
For distinct $i\in B_n$ and $j\in B_m$, separation holds by the choice of $P_{k_n}$ if $n=m$.
If $n<m$, then $|i-j|\geq a_m-b_n\geq a_m-b_{m-1}>Q_m\geq\max(q_{k_n},q_{k_m})$.
The case $m<n$ is symmetric.

\textbf{Case 3: $\cf_{\mathrm{pud}}$.}
Since $\cf_{\mathrm{pud}}$ has the Ramsey property, for each $k$ there exists $r_k\in\{0,\ldots,q_k-1\}$ such that $P_k:=\{n\in E_k\colon n\equiv r_k\pmod{q_k}\}\in\cf_{\mathrm{pud}}$.
Write $\delta_k=\overline d(P_k)>0$.
Recursively choose intervals $I_n=[a_n,b_n]\cap\mathbb N_0$ satisfying $a_n-b_{n-1}>Q_n$ for $n\geq1$ and
\[
|P_{k_n}\cap I_n|
\geq\frac{\delta_{k_n}}4(b_n+1).
\]
Indeed, after choosing $a_n$ sufficiently far to the right, choose $b_n>a_n$ arbitrarily large with $|P_{k_n}\cap[0,b_n]|\geq\delta_{k_n}(b_n+1)/2$ and $a_n\leq\delta_{k_n}(b_n+1)/4$.
Then subtract the at most $a_n$ points in $[0,a_n-1]$.
Also require $b_n\geq n$.
Put $B_n=P_{k_n}\cap I_n$ and $F_k=\bigcup_{\{n\colon k_n=k\}}B_n$.
Along the infinitely many indices with $k_n=k$,
\[
\frac{|F_k\cap[0,b_n]|}{b_n+1}
\geq\frac{|B_n|}{b_n+1}
\geq\frac{\delta_k}{4}.
\]
Since $b_n\to\infty$, we have $\overline d(F_k)\geq\delta_k/4>0$.
Pairwise disjointness and the required separation condition follow by the same argument as in Case 2.

\textbf{Case 4: $\cf_{\mathrm{ps}}$.}
Since $\cf_{\mathrm{ps}}$ has the Ramsey property, for each $k$ there exists $r_k\in\{0,\ldots,q_k-1\}$ such that $P_k:=\{n\in E_k\colon n\equiv r_k\pmod{q_k}\}\in\cf_{\mathrm{ps}}$.
Write $P_k=S_k\cap T_k$, where $S_k$ is syndetic and $T_k$ is thick.
Choose $R_k\in\mathbb N$ such that every interval of $R_k$ consecutive integers meets $S_k$.
Recursively choose finite integer intervals $I_n=[a_n,b_n]\subset T_{k_n}$ such that $b_n-a_n\to\infty$ and $a_n-b_{n-1}>Q_n$ for $n\geq1$.
Put $B_n=P_{k_n}\cap I_n$ and
\[
F_k=\bigcup_{\{n\colon k_n=k\}}B_n.
\]
Fix $k$ and let $n$ satisfy $k_n=k$.
Since $I_n\subset T_k$, we have $B_n=S_k\cap I_n$.
For every $t\in[a_n+R_k,b_n]$, there exists $s\in S_k\cap[t-R_k,t]\subset B_n$, and hence $t\in B_n+[0,R_k]$.
Thus
\[
[a_n+R_k,b_n]
\subset B_n+[0,R_k]
\subset F_k+[0,R_k].
\]
Since infinitely many $n$ satisfy $k_n=k$ and $b_n-a_n\to\infty$, the set $F_k+[0,R_k]$ is thick, so $F_k$ is piecewise syndetic.
The $I_n$ are pairwise disjoint, hence so are the $F_k$.
For distinct $i\in B_n$ and $j\in B_m$, separation holds by the choice of $P_{k_n}$ if $n=m$.
If $n<m$, then $|i-j|\geq a_m-b_n\geq a_m-b_{m-1}>Q_m\geq\max(q_{k_n},q_{k_m})$.
The case $m<n$ is symmetric.
\end{proof}

\section{From ordinary hypercyclicity to weighted hypercyclicity}
In this section, we prove Theorem \ref{thm:shift-criterion} and its simultaneous-weight variants.
The starting point is Shkarin's result on products of hypercyclic operators with compact group rotations.

\begin{proposition}[{\cite[Corollary 4.1]{S08}}]\label{prop:Shkarin}
    Let $T$ be an operator on a complex F-space $X$ and let $G$ be a compact metrizable topological group with a generator $g\in G$.
    If $T$ is hypercyclic, then $\HC (T)\times G=\Trans(T\times R_g)$.
\end{proposition}
\begin{remark}\label{rem:hypercyclic-weak-scattering}
    In Proposition \ref{prop:Shkarin}, since $T$ is hypercyclic, $\HC(T)$ is a dense $G_\delta$ subset of $X$. 
    By Theorem \ref{thm:birkhoff}, $\Trans(T\times R_g)=\HC(T)\times G$ is a dense $G_\delta$ subset of $X\times G$, which gives that $T$ is weakly disjoint from $R_g$.
    That is, $T$ is weakly scattering.
    Thus every hypercyclic operator $T\colon X\to X$ is weakly scattering, where $X$ is a complex F-space.
\end{remark}

For every bounded complex sequence $\bs \lambda$, write
\[
\Sigma_{\bs \lambda}: =\overline{\{B^n\bs \lambda\colon n\in\mathbb N_0\}} 
\quad\text{ and }\quad 
\Omega_{\bs \lambda}: =\bigcap_{N\in\mathbb N_0}\overline{\{B^n\bs \lambda\colon n\geq N\}}.
\]
These are nonempty compact subsets of $\mathbb C^{\mathbb N_0}$ (endowed with the product topology) with $B\Sigma_{\bs\lambda}\subset \Sigma_{\bs\lambda}$ and $B\Omega_{\bs \lambda}=\Omega_{\bs \lambda}$.
The sequence $\bs \lambda$ is a transitive point of $\Sigma_{\bs \lambda}$.
If $\bs \nu\in\Omega_{\bs \lambda}$, then $\Sigma_{\bs \nu}\subset\Omega_{\bs \lambda}$. 
For every compatible metric $d$ on $\mathbb C^{\mathbb N_0}$, one has \[d(B^n\bs\lambda,\Omega_{\bs\lambda})\to 0,\quad n\to\infty, \]
where $d(x,\Omega_{\bs\lambda})=\inf_{y\in\Omega_{\bs\lambda}}d(x,y)$. 
Consequently, for every open subset $O$ of $\Sigma_{\bs\lambda}$ with $\Omega_{\bs\lambda}\subset O$, there exists $N\in\mathbb N_0$ such that $B^n\bs\lambda\in O$ for all $n>N$.
We now apply disjointness to each such orbit closure and then pass from orbit density to return sets in $\cf$.

\begin{proof}[Proof of Theorem \ref{thm:shift-criterion}]
We may assume that $X$ is not a singleton and $T$ is hypercyclic. 
We fix $x\in\HCF(T)$ and a nonempty open subset $U$ of $X$.
For $\bs \nu\in\Omega_{\bs\lambda}$, the system $(\Sigma_{\bs \nu},B)$ is minimal and distal since it is a subsystem of $(\Omega_{\bs\lambda},B)$. 
Moreover, by $\bs 0\notin\Omega_{\bs\lambda}$, we have $\nu_p\neq 0$ for some $p\in\mathbb N_0$.
The continuous map
\[
\Phi\colon X\times\Sigma_{\bs \nu}\to X,
\quad \Phi(u,\bs z)=z_0u,
\]
is onto, since $\Phi(\nu_p^{-1}u,B^p\bs \nu)=u$. 
So $\Phi^{-1}U$ is a nonempty open subset of $X\times \Sigma_{\bs\nu}$. 
By Remarks \ref{rem:hypercyclic-weak-scattering} and \ref{rem:weak-scattering-disjoint-distal}, $T$ is disjoint from $(\Sigma_{\bs \nu},B)$.
By Proposition \ref{prop:compactification-disjoint}, we have that $(x,\bs \nu)$ is transitive for $(X\times \Sigma_{\bs\nu},T\times B)$. 
Hence there is $r_{\bs\nu}\in\mathbb{N}_0$ with 
$(T\times B)^{r_{\bs\nu}}(x,\bs\nu)\in\Phi^{-1}U$. 
That is $\nu_{r_{\bs\nu}}T^{r_{\bs\nu}}x\in U$.

By continuity, there exist a neighborhood $O_{\bs\nu}$ of $\bs\nu$ in $\Sigma_{\bs\lambda}$ and a neighborhood $V_{\bs\nu}$ of $x$ such that
\[
c_{r_{\bs\nu}}T^{r_{\bs\nu}}z\in U
\quad(\bs c\in O_{\bs\nu},\ z\in V_{\bs\nu}).
\]
Now the collection $\{O_{\bs\nu}\colon \nu\in\Omega_{\bs\lambda}\}$ is an open cover of $\Omega_{\bs\lambda}$, which is compact. 
Choose a finite subcover $O_1,\dotsc,O_q$, with corresponding neighborhoods $V_1,\dotsc,V_q$ and integers $r_1,\dotsc,r_q$.
Let $N\in\mathbb{N}_0$ be such that $B^n\bs\lambda\in\bigcup_{j=1}^q O_j$ for every $n>N$.

Since $X$ is not a singleton, the set
\[
V: =\bigcap_{j=1}^qV_j\setminus\{T^lx\colon 0\leq l\leq N\}
\]
is nonempty and open.
Put $A=N_T(x,V)\in\cf$ and $A_j=\{n\in A\colon B^n\bs\lambda\in O_j\}$ for $j=1,2,\dotsc,q$.
Since $A\subset [N+1,\infty)$, we have $A=\bigcup_{j=1}^q A_j$. 
If $n\in A_j$, then
\[
\lambda_{n+r_j}T^{n+r_j}x
=(B^n\bs\lambda)_{r_j}T^{r_j}(T^nx)\in U.
\]
Consequently,
\[
\bigcup_{j=1}^q(r_j+A_j)\subset N_T^{\bs\lambda}(x,U).
\]
By the CuSP property, the left-hand side lies in $\cf$.
Then $N_T^{\bs\lambda}(x,U)\in\cf$, which gives that $x\in\HCF^{\bs\lambda}(T)$. 
\end{proof}

Now we turn to the weighted hypercyclicity of chaotic operators.
An operator is said to be \emph{chaotic} if it is hypercyclic and has a dense set of periodic vectors.
In the preceding proof, distality supplies minimal scalar orbit closures and disjointness from them.
For chaotic operators, the distality assumption can be replaced by minimality.
\begin{theorem}\label{thm:chaotic-uniform-recurrent}
 Let $\bs\lambda$ be a bounded complex sequence.
 If $(\Omega_{\bs\lambda},B)$ is minimal with $\bs 0\notin \Omega_{\bs\lambda}$, then for every complex F-space $X$, every chaotic operator $T\colon X\to X$ and every family $\cf$ with the CuSP property, one has \[\HCF(T)\subset \HCF^{\bs\lambda}(T). \]
\end{theorem}
\begin{proof}
We may assume that $X$ is not a singleton. 
Since $(\Omega_{\bs\lambda},B)$ is minimal, we have $\Sigma_{\bs\nu}=\Omega_{\bs\lambda}$ for every $\bs\nu\in\Omega_{\bs\lambda}$.
Since $T$ is hypercyclic, by Ansari's Theorem (see {\cite[Theorem 1]{A95}}), $T^p$ is hypercyclic for every $p\in\mathbb N$ (i.e.\ $T$ is totally transitive).
Since $T$ has dense periodic vectors, by Remark \ref{rem:totally-transitive-minimal}, it is disjoint from every minimal compact dynamical system.
Hence $T$ is disjoint from $(\Sigma_{\bs\nu},B)$ for every $\bs\nu\in\Omega_{\bs\lambda}$. 
For each $x\in\HCF(T)$ and nonempty open subset $U$ of $X$, this disjointness gives the integers $r_{\bs\nu}$ used in the proof of Theorem \ref{thm:shift-criterion}.
Its neighborhood construction, finite cover and CuSP argument then apply unchanged.
\end{proof}

 Let $\bs\lambda$ be a nonzero bounded complex sequence with $(\Sigma_{\bs\lambda},B)$ minimal.
Then $\Omega_{\bs\lambda}=\Sigma_{\bs\lambda}$, which gives that $(\Omega_{\bs\lambda},B)$ is minimal with $\bs 0\notin \Omega_{\bs\lambda}$.
For example, the following nonconstant weight takes the value zero infinitely often.

\begin{example}
Let $\Sigma_2=\{0,1\}^{\mathbb N_0}$ be the one-sided symbolic space endowed with the product topology, and let $\sigma:\Sigma_2\to\Sigma_2$ be the left shift
defined by
\[
\sigma((x_n)_n)=(x_{n+1})_n.
\]
Consider the Chacon substitution $\chi$ on $\{0,1\}$ given by
\[
\chi(0)=0010,\quad \chi(1)=1.
\]
Let $\bs\lambda\in\Sigma_2$ be the fixed point of $\chi$ beginning with $0$, and set
\[
X_{\bs\lambda}
=\overline{\{\sigma^n\bs\lambda:n\geq0\}},
\]
where the closure is taken in $\Sigma_2$.
The associated subshift $(X_{\bs\lambda},\sigma)$ is minimal and topologically weakly mixing (see \cite[Examples 1.27 and 6.124]{B22}).
Consequently, $(X_{\bs\lambda},\sigma)$ is not distal. 
By Theorems \ref{thm:removal} and \ref{thm:chaotic-uniform-recurrent}, for every complex F-space $X$, every chaotic operator $T\colon X\to X$ and every family $\cf$ with the CuSP property, one has $\HCF(T)=\HCF^{\bs\lambda}(T)$.
\end{example}
There are two different simultaneous-weight questions.
The first concerns a hypercyclic vector for a direct sum, with a separate weight on each coordinate.
Shkarin proved the corresponding ordinary-hypercyclic assertion for geometric unimodular weights \cite[Proposition 4.9]{S08}.
We replace these geometric weights by bounded sequences satisfying a condition on their joint shift system.

Let $\bs\lambda^{(1)},\dotsc,\bs\lambda^{(s)}$ be bounded complex sequences, and put $[\bs\lambda]=\bigl(\bs\lambda^{(1)},\dotsc,\bs\lambda^{(s)}\bigr)$.
Denote $B^{\oplus s}=B\oplus B\oplus \dotsb\oplus B$ ($s$ terms), where $B$ is the backward shift on $\mathbb C^{\mathbb N_0} $.
Denote 
\[
\Sigma_{[\bs \lambda]}: =\overline{\{(B^{\oplus s})^n[\bs \lambda]\colon n\in\mathbb N_0\}} 
\quad\text{ and }\quad 
\Omega_{[\bs \lambda]}: =\bigcap_{N\in\mathbb N_0}\overline{\{(B^{\oplus s})^n[\bs \lambda]\colon n\geq N\}}.
\]
Let $\cf$ be a family.
Let $X_1,\dotsc,X_s$ be a collection of complex F-spaces. 
For every collection of operators $T_1,T_2,\dotsc, T_s$ on $X_1,X_2,\dotsc, X_s$ respectively, put 
\[
    \HCF^{[\bs\lambda]}\biggl(\bigoplus_{i=1}^sT_i\biggr): =\biggl\{x\in \bigoplus_{i=1}^s X_i\colon \bigcap_{i=1}^s N_{T_i}^{\lambda^{(i)}}(x_i,U_i)\in\cf,\ \forall \text{ nonempty open }U_i\subset X_i\biggr\}. 
\]
We write $\HC^{[\bs\lambda]}(\bigoplus_{i=1}^sT_i)$ when $\cf=\cf_{\mathrm{inf}}$.
The next theorem concerns an $\cf$-hypercyclic vector of the direct sum. 
Its scalar hypothesis also requires simultaneous nonvanishing.
For example, the alternating sequences $1,0,1,0,\dotsc$ and $0,1,0,1,\dotsc$ are separately nonzero, but their joint shift orbit has no point whose two zeroth coordinates are both nonzero.

\begin{theorem}\label{thm:multi-sequences}
Let $\bs\lambda^{(1)},\dotsc,\bs\lambda^{(s)}$ be bounded complex sequences, and put $[\bs\lambda]=\bigl(\bs\lambda^{(1)},\dotsc,\bs\lambda^{(s)}\bigr)$.
Suppose that $(\Omega_{[\bs\lambda]},B^{\oplus s})$ is distal and that, for every $[\bs\nu]\in\Omega_{[\bs\lambda]}$, there exists $[\bs \xi]\in\Sigma_{[\bs\nu]}$ such that
\[
\xi_0^{(i)}\neq 0\quad(1\leq i\leq s).
\]
Then for every collection $X_1,X_2,\dotsc, X_s$ of complex F-spaces, every collection of operators $T_1,T_2,\dotsc, T_s$ on $X_1,X_2,\dotsc, X_s$ respectively and every family $\cf$ with the CuSP property, one has
\[
    \HCF\biggl(\bigoplus_{i=1}^sT_i\biggr)\subset \HCF^{[\bs\lambda]}\biggl(\bigoplus_{i=1}^sT_i \biggr).
\]
\end{theorem}
\begin{proof}
Write $X=\bigoplus_{i=1}^s X_i$ and $T=\bigoplus_{i=1}^s T_i$. 
We may assume that $X$ is not a singleton and $T$ is hypercyclic.
Let $x=(x_1,\dotsc,x_s)\in \HCF(T)$ and let
$U_i\subset X_i$, $i=1,\dotsc,s$, be nonempty and open.
Fix $[\bs\nu]=(\bs\nu^{(1)},\dotsc,\bs\nu^{(s)})\in\Omega_{[\bs\lambda]}$.
Since $(\Omega_{[\bs\lambda]},B^{\oplus s})$ is distal,
$(\Sigma_{[\bs\nu]},B^{\oplus s})$ is minimal and distal.

The map 
\[
\Phi\colon X\times\Sigma_{[\bs\nu]}\to X,\quad \Phi(u,[\bs z])=(z_0^{(1)}u_1,\dotsc,z_0^{(s)}u_s), 
\]
is continuous. 
By assumption, there exists
$[\bs \xi]\in\Sigma_{[\bs\nu]}$ such that
$\xi_0^{(i)}\neq 0$ for every $i=1,\dotsc,s$.
For every $u=(u_1,\dotsc,u_s)\in X$, we have
\[
\Phi\left(((\xi_0^{(1)})^{-1}u_1,\dotsc,(\xi_0^{(s)})^{-1}u_s),[\bs \xi]\right)=u,
 \]
which implies that $\Phi$ is onto.
Hence
$\Phi^{-1}(U_1\times\cdots\times U_s)$ is a nonempty open subset of
$X\times\Sigma_{[\bs\nu]}$.

As in the proof of Theorem \ref{thm:shift-criterion}, being weakly scattering,
$T$ is disjoint from the minimal distal system $(\Sigma_{[\bs\nu]},B^{\oplus s})$.
Since $x$ is hypercyclic, $(x,[\bs\nu])$ is transitive for $\left(X\times \Sigma_{[\bs\nu]},T\times B^{\oplus s}\right)$.
Thus there exists
$r_{[\bs\nu]}\in\mathbb N_0$ such that
$(T^{r_{[\bs\nu]}}x,(B^{\oplus s})^{r_{[\bs\nu]}}[\bs\nu])\in
\Phi^{-1}(U_1\times\cdots\times U_s)$.
Equivalently, $\nu_{r_{[\bs\nu]}}^{(i)}T_i^{r_{[\bs\nu]}}x_i\in U_i$ for $i=1,\dotsc,s$.

By continuity, choose neighborhoods $O_{[\bs\nu]}$ of $[\bs\nu]$ in $\Sigma_{[\bs\lambda]}$ and $V_{[\bs\nu]}$ of $x$ such that
\[
c_{r_{[\bs\nu]}}^{(i)}T_i^{r_{[\bs\nu]}}z_i\in U_i\quad (1\leq i\leq s,[\bs c]\in O_{[\bs\nu]},z\in V_{[\bs\nu]}). 
\]
Take a finite subcover $O_1,\dotsc,O_q$ of $\{O_{[\bs\nu]}\colon [\bs\nu]\in\Omega_{[\bs\lambda]}\}$, with corresponding neighborhoods $V_1,V_2,\dotsc, V_q$ and integers $r_1,r_2,\dotsc,r_q$. 
As in the proof of Theorem \ref{thm:shift-criterion}, $(B^{\oplus s})^n[\bs\lambda]\in\bigcup_{j=1}^qO_j$ for all $n>N$, for some $N\in\mathbb N_0$.
Put $V=\bigcap_jV_j\setminus\{T^lx\colon 0\leq l\leq N\}$ and $A=N_T(x,V)\in\cf$.
The sets $A_j=\{n\in A\colon (B^{\oplus s})^n[\bs\lambda]\in O_j\}\ (j=1,2,\dotsc,q)$ cover $A$ and satisfy
\[
\bigcup_{j=1}^q(r_j+A_j)\subset\bigcap_{i=1}^qN_{T_i}^{\bs\lambda^{(i)}}(x_i,U_i).
\]
The conclusion follows from the CuSP property.
\end{proof}

Bayart and Costakis proved a simultaneous-return result for polynomial phases with zero constant terms and frequently hypercyclic vectors \cite[Theorem 3.4]{BC14}.
The second question uses the same vector in every coordinate and asks for simultaneous returns to one open set.
For this purpose we use the stronger sufficient condition that the scalar orbit closure contains a point with a common nonzero zeroth coordinate.

\begin{theorem}\label{thm:multi-sequences-one-operator}
Let $\bs\lambda^{(1)},\dotsc,\bs\lambda^{(s)}$ be bounded complex sequences, and put
$[\bs\lambda]=\bigl(\bs\lambda^{(1)},\dotsc,\bs\lambda^{(s)}\bigr)$.
Suppose that $(\Omega_{[\bs\lambda]},B^{\oplus s})$ is distal and that, for every $[\bs\nu]\in\Omega_{[\bs\lambda]}$, there exists $[\bs \xi]\in\Sigma_{[\bs \nu]}$ such that
\[
\xi_0^{(1)}=\cdots=\xi_0^{(s)}\neq 0.
\]
For every operator $T$ on a complex F-space $X$, every CuSP family $\cf$ and every $x\in\HCF(T)$, one then has
\[
\bigcap_{i=1}^sN_T^{\bs\lambda^{(i)}}(x,U)\in\cf
\]
for every nonempty open subset $U$ of $X$.
\end{theorem}

\begin{proof}
We may assume that $X$ is not a singleton and $T$ is hypercyclic. 
Let $x\in\HCF(T)$ and let $U$ be a nonempty open subset of $X$.
Fix $[\bs\nu]=(\bs\nu^{(1)},\dotsc,\bs\nu^{(s)})\in\Omega_{[\bs\lambda]}$.
Since $(\Omega_{[\bs\lambda]},B^{\oplus s})$ is distal, $(\Sigma_{[\bs\nu]},B^{\oplus s})$ is a minimal distal compact dynamical system.

Consider the continuous map 
\[
    \Phi\colon X\times\Sigma_{[\bs \nu]}\to X^s,\quad \Phi(u,[\bs z])=(z_0^{(1)}u,\dotsc,z_0^{(s)}u). 
\]
By assumption, there exists $[\bs \xi]\in\Sigma_{[\bs\nu]}$ such that
$\xi_0^{(1)}=\dotsb=\xi_0^{(s)}=\colon \xi_0\neq0$.
For every $u\in U$, we have $\Phi(\xi_0^{-1}u,[\bs \xi])=(u,\dotsc,u)\in U^s$.
Hence $\Phi^{-1}(U^s)$ is a nonempty open subset of $X\times\Sigma_{[\bs\nu]}$.

As in the proof of Theorem \ref{thm:shift-criterion}, being weakly scattering, $T$ is disjoint from $(\Sigma_{[\bs\nu]},B^{\oplus s})$.
Since $x$ is hypercyclic, $(x,[\bs\nu])$ is transitive for $(X\times\Sigma_{[\bs\nu]},T\times B^{\oplus s})$.
Thus there exists $r_{[\bs\nu]}\in\mathbb N_0$ such that
$(T^{r_{[\bs\nu]}}x,(B^{\oplus s})^{r_{[\bs\nu]}}[\bs\nu])\in\Phi^{-1}(U^s)$.
Equivalently, $\nu_{r_{[\bs\nu]}}^{(j)}T^{r_{[\bs\nu]}}x\in U$ for every $j=1,\dotsc,s$.

The finite-cover argument in the proof of Theorem \ref{thm:multi-sequences} applies with $V_{[\bs\nu]}\subset X$ and the coordinate conditions
$c_{r_{[\bs\nu]}}^{(i)}T^{r_{[\bs\nu]}}z\in U$ for $i=1,2,\dotsc,s$. 
It produces $A\in\cf$, a finite cover $A=\bigcup_{j=1}^qA_j$ and shifts $r_1,r_2,\dotsc,r_q$ such that
\[
\bigcup_{j=1}^q(r_j+A_j)\subset\bigcap_{i=1}^sN_T^{\bs\lambda^{(i)}}(x,U).
\]
The CuSP property gives the claimed simultaneous return set.
\end{proof}

\section{Applications}
We apply the preceding results by realizing scalar sequences as continuous observations of compact distal systems.
Periodic weights and polynomial phases give elementary instances.
We then pass to nilsequences and turn to simultaneous weights and multiplicative differences.

The following lemma shows that if a bounded complex sequence can be obtained as a continuous observation along the orbit of a distal system, then its backward-shift orbit closure is itself distal.

\begin{lemma}\label{lem:Sigma-distal}
Let $\bs\lambda$ be a bounded complex sequence.
Then the following assertions are equivalent.
\begin{enumerate}
    \item There exist a distal compact dynamical system $(Y,S)$, $y\in Y$ and a continuous function $\varphi\colon Y\to \mathbb C$ such that $\bs\lambda=(\varphi(S^ny))_n$;
    \item $(\Sigma_{\bs\lambda}, B)$ is distal.
\end{enumerate}
\end{lemma}

\begin{proof}
(1)$\Rightarrow$(2). 
Let $Y_0=\overline{\Orb(y,S)}$ and define
\[
\Phi\colon Y_0\to\mathbb C^{\mathbb N_0}, 
\quad \Phi(u)=(\varphi(S^nu))_n. 
\]
This map is continuous, $\Phi(y)=\bs\lambda$ and $\Phi S=B\Phi$.
Since the orbit of $y$ is dense in $Y_0$ and $Y_0$ is compact, its image is exactly $\Sigma_{\bs\lambda}$.
Thus $(\Sigma_{\bs\lambda},B)$ is a factor of the distal system $(Y_0,S)$ and is distal.

(2)$\Rightarrow$(1). 
Consider $(Y,S)=(\Sigma_{\bs\lambda},B)$, $y=\bs\lambda$ and the continuous map $\varphi(\bs z)=z_0$.
\end{proof}

\begin{remark}\label{rem:Sigma-closed}
Let $\bs\lambda$ and $\bs\mu$ be two bounded complex sequences such that both $(\Sigma_{\bs\lambda},B)$ and $(\Sigma_{\bs\mu},B)$ are distal.
Then for every continuous function $\Phi\colon (\mathbb C^{\mathbb N_0})^2\to\mathbb C^{\mathbb N_0}$ with $\Phi (B\times B)=B \Phi$, $(\Sigma_{\Phi(\bs\lambda,\bs\mu)},B)$ is distal.
Indeed, let $p_0$ be the zeroth-coordinate projection on $\mathbb C^{\mathbb N_0}$ and define $\varphi=p_0\circ\Phi$ on $\Sigma_{\bs\lambda}\times\Sigma_{\bs\mu}$.
Since $\Phi$ is continuous and commutes with the backward shift, we have
\[
\Phi(\bs\lambda,\bs\mu)=\bigl(\varphi((B\times B)^n(\bs\lambda,\bs\mu))\bigr)_n.
\]
Since the product $(\Sigma_{\bs\lambda}\times\Sigma_{\bs\mu},
B\times B)$ is distal and $\varphi$ is continuous, the conclusion follows from Lemma \ref{lem:Sigma-distal}.
In particular, the conclusion applies to sums $\Phi(\bs\lambda,\bs\mu)=\bs\lambda+\bs\mu$, products $\Phi(\bs\lambda,\bs\mu)=\bs\lambda\bs\mu: =(\lambda_n\mu_n)_n$, complex conjugation $\Phi(\bs\lambda,\bs\mu)=\overline{\bs\lambda}: =(\overline{\lambda_n})_n$, and more generally to every complex polynomial $\Phi(\bs\lambda,\bs\mu)=\sum_{p,q}a_{p,q}\bs\lambda^p\bs\mu^q$.
\end{remark}

\begin{remark}\label{rem:shift-criterion-sigma}
If $\bs\lambda\ne\bs0$ and $(\Sigma_{\bs\lambda},B)$ is distal, then
$\Omega_{\bs\lambda}=\Sigma_{\bs\lambda}$ by minimality.
This set cannot contain the fixed point $\bs0$, since otherwise minimality would force it to be $\{\bs0\}$.
Thus Theorem \ref{thm:shift-criterion} applies.
Distality of the full orbit closure is stronger than the condition on the $\omega$-limit set.
For example, if $\mu_n=(n+2)/(n+1)$, then $B^n\bs\mu\to\bs1$ as $n\to\infty$ in the product topology, where $\bs1=(1,1,\dotsc)$.
Hence $\Omega_{\bs\mu}=\{\bs1\}$ is distal, while $\inf_{n\in\mathbb N_0}d(B^n\bs\mu,B^n\bs 1)=0$ for every compatible metric $d$ on $\mathbb C^{\mathbb N_0}$, giving that $\Sigma_{\bs\mu}$ is not distal. 
\end{remark}

The next examples recover powers and unimodular rotations.
These are the assertions of \cite[Theorem 8.8]{BGLP22}, whose right-invariance hypothesis is already implied by CuSP.
The powers assertion also uses the stated dilation condition.

\begin{example}\label{exam:periodic}\
Let $p\in \mathbb{N}$ and let $\bs\lambda$ be given by $\lambda_n=1$ if $p|n$, and $\lambda_n=0$ otherwise.
The shift permutes the finitely many translates of this periodic sequence, so $(\Sigma_{\bs\lambda},B)$ is distal.
By Remark \ref{rem:shift-criterion-sigma}, $(\Omega_{\bs\lambda},B)$ is distal with $\bs 0\notin \Omega_{\bs\lambda}$.
By Theorems \ref{thm:removal} and \ref{thm:shift-criterion}, for every complex F-space $X$, every operator $T\colon X\to X$ and every family $\cf$ with the CuSP property, we have $\HCF(T)=\HCF^{\bs\lambda}(T)$.
Assume in addition that for every $t\in\mathbb N$ and every $A\subset \mathbb N_0$, $tA\in\cf$ if and only if $A\in\cf$. 
For every open subset $U$ of $X$ with $0\notin U$, we have $N_T^{\bs\lambda}(x,U)=pN_{T^p}(x,U)$.
It follows that $\HCF^{\bs\lambda}(T)=\HCF(T^p)$.
Thus $\HCF(T)=\HCF(T^p)$. 
\end{example}
 
\begin{example}\label{exam:rotation}
    Denote $\mathbb T=\{z\in\mathbb C\colon |z|=1\}$ and let $\lambda\in\mathbb T$.
    Consider the map $R\colon \T\to\T$ given by \[R z=\lambda z, \quad z\in\T. \]
Then $(\T,R)$ is distal and $\lambda^n=R^n1$ for every $n\in\mathbb N_0$.
By Lemma \ref{lem:Sigma-distal}, $(\Sigma_{(\lambda^n)_n},B)$ is distal.
Remark \ref{rem:shift-criterion-sigma} and Theorems \ref{thm:removal} and \ref{thm:shift-criterion} give, for every operator on a complex F-space and every CuSP family,
$\HCF(T)=\HCF^{(\lambda^n)_n}(T)$.
That is, $\HCF(T)=\HCF(\lambda T)$.
\end{example}
\begin{remark}
    In Example \ref{exam:periodic}, if for every $t\in\mathbb N$ and every $A\subset \mathbb N_0$, $tA\in\cf$ if and only if $A\in\cf$, then we obtain $\HCF(T)=\HCF(T^p)$.
    Let $\bs\mu$ be a complex sequence whose nonzero values are bounded and bounded away from zero such that $(\Omega_{\bs\mu},B)$ is distal with $\bs 0\notin \Omega_{\bs\mu}$.
    By Theorems \ref{thm:removal} and \ref{thm:shift-criterion}, we have $\HCF(T)=\HCF^{\bs\mu}(T)$ and $\HCF(T^p)=\HCF^{\bs\mu}(T^p)$.
    Hence we obtain $\HCF^{\bs\mu}(T)=\HCF^{\bs\mu}(T^p)$.
    In the same vein, following Example \ref{exam:rotation}, we have $\HCF^{\bs\mu}(T)=\HCF^{\bs\mu}(\lambda T)$.
\end{remark}
Polynomial phases of higher degree are obtained from affine skew rotations.
They can be realized as continuous observations of distal affine systems on $(\mathbb R/\mathbb Z)^d$.
This allows Theorems \ref{thm:removal} and \ref{thm:shift-criterion} to recover and extend Theorem \ref{thm:bc14} of Bayart and Costakis from ordinary hypercyclicity to families with the CuSP property.
For every $x\in\mathbb R$, let $|x|_{\mathbb R/\mathbb Z}=\min\{|x-k|\colon k\in\mathbb Z\}$.
For every real polynomial $P$, we denote $\bs e_P=(e^{2\pi iP(n)})_n$, the unimodular sequence with $P$ as its phase.
\begin{example}\label{exam:polynomial}
Let $P$ be a real polynomial.
If $P$ is constant, then $\bs e_P$ is constant and $(\Sigma_{\bs e_P},B)$ is distal.
Otherwise let $d=\deg P\geq1$ and write $P(n)=\sum_{i=0}^da_i\binom ni$.
Consider $Y=(\mathbb R/\mathbb Z)^d$ and define the transformation $S\colon Y\to Y$ by
\[
S(x_1,\dotsc,x_d)=(x_1+a_d,x_2+x_1,x_3+x_2,\dotsc,x_d+x_{d-1})\pmod 1
\]
for every $(x_1,\dotsc,x_d)\in Y$.
Considering the continuous function $\varphi\colon Y\to\mathbb C$ given by \[\varphi (x_1, x_2,\dotsc, x_d)=e^{2\pi i x_d},\quad (x_1,\dotsc,x_d)\in Y, \]
and the point $y\in Y$ defined by $y_j=a_{d-j}$ for $1\leq j\leq d$, the binomial identity gives \[\varphi(S^n y)=e^{2\pi iP(n)},\quad n\in\mathbb N_0. \]
The system $(Y,S)$ and the corresponding orbit computation leading to the polynomial $P$ are discussed by Furstenberg, see the discussion preceding \cite[Theorem 1.6]{F81}.

Consider the metric $d_Y$ on $Y$ given by $d_Y(u,v)=\max \{|u_j-v_j|_{\mathbb R/\mathbb Z}\colon 1\leq j\leq d\}$ for $u,v\in Y$.
Let $u,v\in Y$ be such that $u\neq v$.
Put $r=\min\{1\leq j\leq d\colon u_j\neq v_j\}$.
Then for every $j<r$, $u_j=v_j$, which gives that $(S^n u)_j=(S^n v)_j$ for every $n\in\mathbb N_0$.
Therefore \[d_Y(S^n u,S^n v)\geq |(S^n u-S^n v)_r|_{\mathbb R/\mathbb Z}=|u_r-v_r|_{\mathbb R/\mathbb Z}\neq 0,\quad n\in\mathbb N_0. \]
Hence, $(Y,S)$ is distal.
By Lemma \ref{lem:Sigma-distal}, $(\Sigma_{\bs e_P},B)$ is distal.
Thus $(\Sigma_{\bs e_P},B)$ is distal for every real polynomial $P$.
By Remark \ref{rem:shift-criterion-sigma}, $(\Omega_{\bs e_P},B)$ is distal with $\bs 0\notin \Omega_{\bs e_P}$.
By Theorems \ref{thm:removal} and \ref{thm:shift-criterion}, for every complex F-space $X$, every operator $T\colon X\to X$ and every family $\cf$ with the CuSP property, we have $\HCF(T)=\HCF^{\bs e_P}(T)$.
\end{example}

In fact, the periodic $0$-$1$-sequence in Example \ref{exam:periodic}, the unimodular geometric sequence $(\lambda^n)_n$ in Example \ref{exam:rotation}, and the polynomial phase sequence $\bs e_P$ in Example \ref{exam:polynomial} can all be placed within a common framework, namely that of nilsequences.
Let $G$ be a nilpotent Lie group and let $\Lambda$ be a discrete cocompact subgroup of $G$.
Consider the compact nilmanifold $G/\Lambda$ and fix $g\in G$.
Let $T_g\colon G/\Lambda\to G/\Lambda$ be given by
\[
  T_g(v\Lambda)=gv\Lambda, \quad v\Lambda\in G/\Lambda.
\]
The compact dynamical system $(G/\Lambda,T_g)$ is called a \emph{nilsystem}.
Every nilsystem is distal, see e.g.\ \cite[2.3]{HKM10}.
Let $u\in G/\Lambda$, and let $\varphi\colon G/\Lambda\to\mathbb C$ be a continuous function.
The sequence $(\varphi(T_g^n u))_n$ is called a \emph{basic nilsequence}.
Every basic nilsequence lies in $\ell^\infty$, the space of bounded complex sequences with norm $\|\bs a\|_\infty=\sup_{n\in\mathbb{N}_0}|a_n|$.
A \emph{nilsequence} is a uniform limit of basic nilsequences.
See \cite{BHK05} and \cite{HKM10} for more details of nilsequences.

\begin{theorem}\label{thm:nilsequence}
Let $\bs\lambda$ be a complex nilsequence which is not identically zero.
 Then for every complex F-space $X$, every operator $T\colon X\to X$ and every family $\cf$ with the CuSP property, we have
 \[
 \HCF(T)\subset \HCF^{\bs\lambda} (T).
 \]
Moreover, if the nonzero values of $\bs\lambda$ are bounded away from zero, then
\[
\HCF(T)=\HCF^{\bs\lambda} (T).
\]
\end{theorem}

\begin{proof}
Choose basic nilsequences $(\bs\lambda^{(j)})_j$ with
$\|\bs\lambda^{(j)}-\bs\lambda\|_\infty\to0$ as $j\to\infty$.
Since every nilsystem is distal, by Lemma \ref{lem:Sigma-distal}, every system $(\Sigma_{\bs\lambda^{(j)}},B)$ is distal. 
Put 
\[\Sigma=\prod_{j=0}^\infty \Sigma_{\bs\lambda^{(j)}}, \quad  \mathcal B=\prod_{j=0}^\infty B, \quad [\bs\lambda]=(\bs\lambda^{(j)})_j. \]
Set $Z=\overline{\Orb([\bs\lambda],\mathcal B)}$.
Then $(Z,\mathcal B)$ is a distal compact dynamical system.
The functions
\[
\Phi_j\colon Z\to\mathbb C,
\quad \Phi_j([\bs z])=z_0^{(j)},
\]
are continuous.
Since the orbit of $[\bs\lambda]$ is dense in $Z$,
\[
\|\Phi_j-\Phi_k\|_{C(Z)}
=\sup_{n\in\mathbb{N}_0}|\lambda_n^{(j)}-\lambda_n^{(k)}|
=\|\bs\lambda^{(j)}-\bs\lambda^{(k)}\|_\infty.
\]
They therefore converge uniformly to a continuous function $\Phi$ on $Z$.
For each $n\in\mathbb{N}_0$, $\Phi(\mathcal B^n[\bs\lambda])=\lambda_n$.
Lemma \ref{lem:Sigma-distal} shows that $(\Sigma_{\bs\lambda},B)$ is distal.
Since $\bs\lambda\ne\bs0$, Remark \ref{rem:shift-criterion-sigma} and Theorem \ref{thm:shift-criterion} yield the first inclusion.
The reverse inclusion under the lower bound follows from Theorem \ref{thm:removal}.
\end{proof}
For instance, if $\alpha$ is irrational, $\lambda_n=1+e^{2\pi in\alpha}$ is a basic nilsequence with no zero terms but with $\inf_{n\in\mathbb N_0}|\lambda_n|=0$.
The first inclusion still holds, whereas Theorem \ref{thm:removal} shows that the reverse inclusion cannot hold for all operators on complex F-spaces. 
The proof of Theorem \ref{thm:nilsequence} also shows that the class of bounded complex sequences with distal orbit closures is closed under uniform convergence. 

We next give an explicit basic nilsequence on the Heisenberg nilmanifold to illustrate Theorem \ref{thm:nilsequence}.
The observation below has modulus between $1$ and $3$.
For every $t\in\mathbb R$, denote $\lfloor t\rfloor=\max\{n\in\mathbb Z\colon n\leq t\}$ and $e(t)=e^{2\pi it}$.

\begin{example}
Let
\[
G=\{[x,y,z]\colon (x,y,z)\in\mathbb R^3\}
\]
be the real Heisenberg group equipped with the
multiplication law
\[
[x,y,z][a,b,c]=[x+a,y+b,z+c+xb].
\]
Then $G$ is a connected, simply connected,
two-step nilpotent Lie group.
The subgroup
\[
\Lambda=\{[x,y,z]\colon (x,y,z)\in\mathbb Z^3\}
\]
is a cocompact lattice in $G$, and hence
$G/\Lambda$ is a compact nilmanifold
(see e.g.\ \cite[Section 4.2.2]{BHK05}).

Define $\varphi\colon G/\Lambda\to\mathbb C$ by
\[
\varphi([x,y,z]\Lambda)
=\sin(2\pi y)e(z-\lfloor y\rfloor x)+2,
\]
where $e(t)=e^{2\pi it}$.
It is not hard to see that $\varphi$ is well-defined and continuous.
For $\alpha,\beta\in \mathbb R$, let $g=[2\beta,\alpha,0]$ and let $u=[0,0,0]\Lambda$.
Since $g^n=[2n\beta,n\alpha,n(n-1)\alpha\beta]$, we have for every $n\in\mathbb{N}_0$
\[
\varphi(T_g^nu)=\sin(2\pi n\alpha)e(n(n-1)\alpha\beta-2n\beta \lfloor n\alpha \rfloor )+2.
\]
For each $n\in\mathbb N_0$, put $\lambda_n=\varphi(T_g^nu)$.
Then we have that $\bs\lambda=(\lambda_n)_n$ is a basic nilsequence which is bounded away from zero.
By Theorem \ref{thm:nilsequence}, for every complex F-space $X$, every operator $T\colon X\to X$ and every family $\cf$ with the CuSP property, we have $\HCF(T)=\HCF^{\bs\lambda}(T)$.
\end{example}

Returning to simultaneous weights, polynomial phases satisfy the joint hypotheses in Section 4.
\begin{example}\label{exam:direct-sum}
    Let $P_1, P_2,\dotsc, P_s$ be real polynomials.
    We have shown in Example \ref{exam:polynomial} that $(\Sigma_{\bs e_{P_i}},B)$ is distal for $i=1,2,\dotsc,s$.
    If we denote
    \[[\bs\lambda]=(\bs e_{P_1},\bs e_{P_2},\dotsc,\bs e_{P_s}), \]
    then $(\Sigma_{[\bs\lambda]}, B^{\oplus s})$ is minimal and distal since it is a subsystem of the product $\prod_{i=1}^s(\Sigma_{\bs e_{P_i}},B)$.
    So $\Sigma_{[\bs\lambda]}=\Omega_{[\bs\lambda]}$. 
    Since $\bs e_{P_i}$ is a unimodular sequence for $i=1,2,\dotsc,s$, by Theorem \ref{thm:multi-sequences}, for every collection $X_1,X_2,\dotsc, X_s$ of complex F-spaces, every collection of operators $T_1,T_2,\dotsc, T_s$ on $X_1,X_2,\dotsc, X_s$ respectively and every family $\cf$ with the CuSP property, one has 
\[
    \HCF\biggl(\bigoplus_{i=1}^sT_i\biggr)\subset \HCF^{[\bs\lambda]}\biggl(\bigoplus_{i=1}^sT_i\biggr).
\]
\end{example}
A related simultaneous-return statement follows from Theorem \ref{thm:multi-sequences-one-operator}.
For a single vector, the zero constant terms provide the common scalar value required by Theorem \ref{thm:multi-sequences-one-operator}.
This yields the following extension of \cite[Theorem 3.4]{BC14} from positive lower density to CuSP families.
\begin{example}
Let $P_1,P_2,\dotsc, P_s$ be real polynomials with zero constant terms. 
Put
\[
[\bs \lambda]=(\bs e_{P_1},\dotsc,\bs e_{P_s}).
\]
In Example \ref{exam:direct-sum}, we obtain that $(\Sigma_{[\bs\lambda]}, B^{\oplus s})$ is minimal and distal with $\Sigma_{[\bs\lambda]}=\Omega_{[\bs\lambda]}$. 
Then for every $[\bs \nu]\in\Omega_{[\bs\lambda]}$, $[\bs \lambda]\in\Sigma_{[\bs\lambda]}=\Omega_{[\bs \lambda]}=\Sigma_{[\bs \nu]}$, whose zeroth coordinates are all $1$ since $P_i(0)=0$.
By Theorem \ref{thm:multi-sequences-one-operator}, for every complex F-space $X$, every operator $T\colon X\to X$, every family $\cf$ with the CuSP property and every $x\in\HCF(T)$,
\[
\bigcap_{i=1}^s N_T^{\bs e_{P_i}}(x,U)\in\cf
\]
for every nonempty open subset $U$ of $X$.
\end{example}

Define the continuous multiplicative-difference map
\[
D\colon\mathbb T^{\mathbb N_0}\to\mathbb T^{\mathbb N_0},
\quad (D\bs\mu)_n=\mu_{n+1}/\mu_n.
\]
It commutes with $B$.
For a real sequence $f$, write $\Delta f(n)=f(n+1)-f(n)$.
Then
\[
D\bigl((e^{2\pi if(n)})_n\bigr)=(e^{2\pi i\Delta f(n)})_n.
\]
In particular, $D^{d+1}\bs e_P=\bs1$ whenever $P$ has degree at most $d$, including constant and zero polynomials.
We shall allow the iterated difference to tend to $1$ rather than equal it.
The following lemma relates this condition to the phase-increment hypothesis of Bayart and Costakis \cite[Theorem 4.3]{BC14}.
Their weights are written as $e^{if(n)}$, and we use the normalization $e^{2\pi if(n)}$.
The degree $d$ in their increment condition corresponds to $r=d+1$ below.

\begin{lemma}\label{lem:difference-increments}
Let $\bs\mu\in\mathbb T^{\mathbb N_0}$ and let $r\geq1$ be an integer.
The following assertions are equivalent.
\begin{enumerate}
\item $(D^r\bs\mu)_n\to1$ as $n\to\infty$; 
\item there is a real sequence $f$ such that $\mu_n=e^{2\pi if(n)}$ for every $n\in\mathbb{N}_0$ and $\Delta^rf(n)\to0$ as $n\to\infty$; 
\item there exist a real sequence $f$ with $\mu_n=e^{2\pi if(n)}$ for every $n\in\mathbb{N}_0$ and real sequences $g_0,\dotsc,g_{r-1}$ such that
\[
f(n+k)-f(n)=\sum_{\ell=0}^{r-1}g_\ell(n)k^\ell+\varepsilon_k(n),
\quad n\in\mathbb{N}_0,\ k\geq1,
\]
where $\varepsilon_k(n)\to0$ as $n\to\infty$ for every fixed positive integer $k$.
\end{enumerate}
Moreover, for any fixed real sequence $f$, the increment condition in {\rm(3)} is equivalent to $\Delta^rf(n)\to0$ as $n\to\infty$.
\end{lemma}
\begin{proof}
(1)$\Rightarrow$(2). 
Choose real $\eta_n\to0$ as $n\to\infty$ such that $e^{2\pi i\eta_n}=(D^r\bs\mu)_n$.
Choose real numbers $f(j)$ with $e^{2\pi if(j)}=\mu_j$ for $0\leq j<r$, and define the remaining terms by $\Delta^rf(n)=\eta_n$.
This is possible since the coefficient of $f(n+r)$ in $\Delta^rf(n)$ is $1$.
The sequence $\widetilde\mu_n=e^{2\pi if(n)}$ has the same first $r$ terms as $\bs\mu$ and satisfies $D^r\bs{\widetilde\mu}=D^r\bs\mu$.
The formula
\[
(D^r\bs\mu)_n=\prod_{j=0}^r\mu_{n+j}^{(-1)^{r-j}\binom rj}
\]
determines $\mu_{n+r}$ from the preceding $r$ terms and $(D^r\bs\mu)_n$.
Induction therefore gives $\widetilde\mu_n=\mu_n$ for every $n$, proving (2).

(2)$\Rightarrow$(1). This follows from $(D^r\bs\mu)_n=e^{2\pi i\Delta^rf(n)}$.

(2)$\Rightarrow$(3). Let $f$ be a real sequence with $\Delta^rf(n)\to0$ as $n\to\infty$.
For each fixed positive integer $k$, the discrete Taylor formula gives
\[
f(n+k)-f(n)=\sum_{j=1}^{r-1}\binom kj\Delta^jf(n)+R_k(n).
\]
Here $R_k(n)=0$ if $k<r$, while for $k\geq r$,
\[
R_k(n)=\sum_{\ell=0}^{k-r} \binom{k-1-\ell}{r-1}\Delta^rf(n+\ell).
\]
Thus $R_k(n)\to0$ as $n\to\infty$ for every fixed $k$.
Since $\binom kj$ is a polynomial in $k$ of degree $j$, this proves the increment condition in (3).
When $r=1$, the first sum is empty and we take $g_0=0$.

(3)$\Rightarrow$(2). Suppose a real sequence $f$ satisfies that increment condition.
For fixed $n$, the polynomial $\sum_{\ell=0}^{r-1}g_\ell(n)k^\ell$ has degree at most $r-1$ in $k$.
Taking its $r$th difference at $k=1$, one has
\[
\Delta^rf(n+1) =\sum_{j=0}^r(-1)^{r-j}\binom rj\varepsilon_{j+1}(n) \to 0,\quad n\to\infty. \qedhere 
\]
\end{proof}

Thus the multiplicative-difference and phase-increment conditions describe the same class of unimodular weights.
The choice of the phase in assertions (2) and (3) in Lemma~\ref{lem:difference-increments} is essential, since adding an arbitrary integer-valued sequence to $f$ leaves $\bs\mu$ unchanged but may destroy either condition.
We now combine these unimodular weights with the weights having distal $\omega$-limit systems considered above.

\begin{theorem}\label{thm:Dr-product}
Let $\bs\lambda$ be a bounded complex sequence such that $(\Omega_{\bs\lambda},B)$ is distal and $\bs0\notin\Omega_{\bs\lambda}$.
Let $\bs\mu$ be unimodular and suppose that $(D^r\bs\mu)_n\to1$ as $n\to\infty$ for some $r\in\mathbb N$, or equivalently that $\bs\mu$ has a real phase satisfying the increment condition in Lemma \ref{lem:difference-increments}.
Then for every complex F-space $X$, every operator $T\colon X\to X$ and every family $\cf$ with the CuSP property, one has \[
\HCF(T)\subset \HCF^{\bs\lambda\bs\mu}(T). 
\]
Moreover, if the nonzero values of $\bs\lambda$ are bounded away from zero, then 
\[
\HCF(T)=\HCF^{\bs\lambda\bs\mu}(T). 
\]
\end{theorem}
\begin{proof}
For $\bs\nu\in\Omega_{\bs\mu}$, choose $n_j\to\infty$ with $B^{n_j}\bs\mu\to\bs\nu$.
Continuity and commutation give
\[
D^r\bs\nu=\lim_{j\to\infty}B^{n_j}D^r\bs\mu=\bs 1.
\]
Put $a_j=(D^j\bs\nu)_0$ for $0\leq j<r$.
Since $D^r\bs\nu=\bs1$, repeated use of the identity
$(D^j\bs\nu)_{n+1}=(D^j\bs\nu)_n(D^{j+1}\bs\nu)_n$ yields
\[
\nu_n=\prod_{j=0}^{r-1}a_j^{\binom nj}.
\]
Choosing real $\theta_j$ with $a_j=e^{2\pi i\theta_j}$ shows that $\bs\nu=\bs e_P$ for
$P(t)=\sum_{j=0}^{r-1}\theta_j\binom tj$.
By Example \ref{exam:polynomial}, $\Sigma_{\bs\nu}$ is minimal and distal.

It follows that $\Omega_{\bs\mu}$ is distal.
Indeed, otherwise two distinct points $\bs b,\bs c\in\Omega_{\bs\mu}$ would satisfy that $\inf_{n\in\mathbb N_0}d(B^n\bs b,B^n\bs c)=0$.
By compactness, along a subsequence their iterates would converge to a common point $\bs z$.
Then $\bs z\in\Sigma_{\bs b}\cap\Sigma_{\bs c}$, and minimality gives
$\Sigma_{\bs b}=\Sigma_{\bs c}$. 
This gives that $(\Sigma_{\bs b},B)$ is not distal, a contradiction. 

Coordinatewise multiplication defines a continuous equivariant map
\[
\Phi\colon\Omega_{\bs\lambda}\times\Omega_{\bs\mu}\to\mathbb C^{\mathbb N_0},
\quad \Phi(\bs a,\bs b)=\bs a\bs b.
\]
It can be checked that $\Omega_{\bs\lambda\bs\mu}\subset \Phi(\Omega_{\bs\lambda}\times\Omega_{\bs\mu})$. 
Thus $\Omega_{\bs\lambda\bs\mu}$ is distal since it is a subsystem of $(\Phi(\Omega_{\bs\lambda}\times\Omega_{\bs\mu}), B)$, which is a factor of $(\Omega_{\bs\lambda}\times\Omega_{\bs\mu},B\times B)$. 
Since all coordinates of every $\bs b\in\Omega_{\bs\mu}$ have modulus one, $\bs 0\notin \Phi(\Omega_{\bs\lambda}\times\Omega_{\bs\mu})$. 
Theorem \ref{thm:shift-criterion} gives the first inclusion.
Finally, multiplication by $\bs\mu$ preserves the moduli of the weights, so Theorem \ref{thm:removal} gives equality under the lower-bound assumption.
\end{proof}

For $\bs\lambda=\bs1$ and ordinary hypercyclicity on Banach spaces, Lemma \ref{lem:difference-increments} identifies Theorem \ref{thm:Dr-product} with \cite[Theorem 4.3]{BC14}.
Theorem \ref{thm:Dr-product} applies to all CuSP families on complex F-spaces and permits multiplication by every bounded complex weight satisfying the distal $\omega$-limit hypothesis, without a compatibility assumption between the two weights.
The following example includes the slowly varying phases considered in \cite{BC14}.

\begin{example}
Let $\alpha,\beta\in\mathbb R$ and put
\[
\mu_n=e^{2\pi i(n+1)^\alpha(\log(n+2))^\beta}.
\]
Choose an integer $r>\max\{\alpha,0\}$.
For $f(t)=(t+1)^\alpha(\log(t+2))^\beta$, differentiation gives
$f^{(r)}(t)=O(t^{\alpha-r}(\log t)^\beta)$ as $t\to\infty$, and hence $f^{(r)}(t)\to0$ as $t\to\infty$.
Using
\[
\Delta^rf(n)=\int_{[0,1]^r}f^{(r)}(n+t_1+\cdots+t_r)\,dt_1\cdots dt_r,
\]
we obtain $(D^r\bs\mu)_n\to1$ as $n\to\infty$.
Therefore Theorem \ref{thm:Dr-product} applies after multiplication by every bounded $\bs\lambda$ whose $\omega$-limit system is distal and avoids $\bs0$.
It gives $\HCF(T)\subset\HCF^{\bs\lambda\bs\mu}(T)$ for every complex F-space $X$, every operator $T\colon X\to X$ and every family $\cf$ with the CuSP property, with equality if the nonzero values of $\bs\lambda$ are bounded away from zero.
\end{example}

\medskip 
\noindent \textbf{AI Declaration}:
This work grew out of the second author's Ph.D. research project, which began in September 2025 under the supervision of the first and third authors.
The authors used AI tools solely to assist with language editing. 
No AI tools were used to derive mathematical results or construct proofs. 
All mathematical results and proofs presented in this paper were developed by the authors, 
who take full responsibility for the mathematical content.

\medskip 
\noindent \textbf{Acknowledgments}: 
The authors were partially supported by the National Key R\&D Program of China (2024YFA1013601), National Natural Science Foundation of China (12222110) and a grant from the Guangdong Provincial Department of Education (2025KCXTD013).

\end{document}